\documentclass[amsfonts, 12pt]{article}

\usepackage{amssymb,amsmath}

\newtheorem{theorem}{Theorem}[section]
\newtheorem{corollary}[theorem]{Corollary}
\newtheorem{example}[theorem]{Example}

\newtheorem{lemma}[theorem]{Lemma}

\newtheorem{remark}[theorem]{Remark}

\def\bR{\mathbb{R}}

\def\bC{\mathbb{C}}

\def\cE{\mathcal{E}}
\def\cF{\mathcal{F}}

\def\cH{\mathcal{H}}
\def\cL{\mathcal{L}}

\def\cP{\mathcal{P}}

\def\cS{\mathcal{S}}

\def\cL{\mathcal{L}}
\def\cZ{\mathcal{Z}}

\begin{document}

\title{Some linear SPDEs 
driven by a fractional noise with Hurst index greater than 1/2}

\author{Raluca M. Balan\footnote{Department of Mathematics and Statistics, University of Ottawa,
585 King Edward Avenue, Ottawa, ON, K1N 6N5, Canada. E-mail
address: rbalan@uottawa.ca} \footnote{Research supported by a
grant from the Natural Sciences and Engineering Research Council
of Canada.}}

\date{February 18, 2011}

\maketitle

\begin{abstract}
\noindent In this article, we identify the necessary and
sufficient conditions for the existence of a random field solution
for some linear s.p.d.e.'s of parabolic and hyperbolic type. These
equations rely on a spatial operator $\cL$ given by the
$L^2$-generator of a $d$-dimensional L\'evy process $X=(X_t)_{t
\geq 0}$, and are driven by a spatially-homogeneous Gaussian
noise, which is fractional in time with Hurst index $H>1/2$. As an application, we consider the case
when $X$ is a $\beta$-stable process, with $\beta \in (0,2]$. In
the parabolic case, we develop a connection with the potential
theory of the Markov process $\bar{X}$ (defined as the
symmetrization of $X$), and we show that the existence of the
solution is related to the existence of a ``weighted''
intersection local time of two independent copies of $\bar{X}$.
\end{abstract}



\vspace{5mm}


\newpage

\section{Introduction}

In 1944, in his seminal article \cite{ito44}, It\^o introduced the
stochastic integral with respect to the Brownian motion, which
turned out to be one of the most fruitful ideas in mathematics in
the 20th century. This lead to the theory of
diffusions (whose origins can be traced back to \cite{ito51}), and
the development of the stochastic calculus with respect to
martingales (initiated in \cite{kunita-watanabe67}). These ideas
have grown into a solid branch of probability theory
called stochastic analysis, which includes 
the study of stochastic partial differential equations (s.p.d.e.'s)

Traditionally, there have been several approaches for the study of
s.p.d.e.'s. The most important are: the Walsh approach which
relies on stochastic integrals with respect to martingale-measures
(see \cite{walsh86}), the Da Prato and Zabczyk approach which uses
stochastic integrals with respect to Hilbert-space-valued Wiener
processes (see \cite{daprato-zabczyk92}), and the Krylov approach
which uses the concept of function-space-valued solution (see
\cite{krylov99}). These approaches have been developed at the same
time, and nowadays a lot of effort is dedicated to unify them (see the
recent survey \cite{dalang-quer10} and the references therein).

The fractional Brownian motion (fBm) was introduced by Kolmogorov in \cite{kolmogorov40}, who called it the ``Wiener spiral'', and is defined as a zero-mean Gaussian process $(B_t)_{t \geq 0}$ with covariance:
$$R_{H}(t,s)=E(B_t B_s)=\frac{1}{2}(t^{2H}+s^{2H}-|t-s|^{2H}).$$
The parameter $H$ lies in $(0,1)$, and is called the Hurst index (due to \cite{hurst51}). The case $H=1/2$ corresponds to the Brownian motion, whereas the cases $H>1/2$ and $H<1/2$ have many contrasting properties and cannot be handled simultaneously. The representation of the fBm as a stochastic integral with respect to the Brownian motion on $\bR$ was obtained as early as 1968 (see \cite{mandelbrot-vanness68}), but the fBm began to be used
intensively in stochastic analysis only in the late 1990's. It is the flexibility which stems from the choice of the parameter $H$ that makes the fBm a much more attractive model for the noise than the Brownian motion (and its infinite-dimensional counterparts).

Among the fBm's remarkable properties is the fact that it is {\em
not} a semimartingale. Consequently, It\^o calculus cannot be used
in this case. A stochastic calculus with respect to the fBm was
developed for the first time in \cite{decreusefond-ustunel98}.
Subsequent important contributions were made in
\cite{alos-mazet-nualart01}, \cite{carmona-coutin-montseny03} and
\cite{duncan-hu-pasik00}. The stochastic integral used by these
authors is an extension of the It\^{o} integral introduced by
Hitsuda in \cite{hitsuda72} (and refined in \cite{kabanov75} and
\cite{skorohod75}), and coincides with the divergence operator.
These techniques are based on Malliavin calculus. Alternative
methods for defining a stochastic integral with respect to the fBm
exploit the H\"older continuity property of its sample paths and
are based on generalized Stieltjes integrals. We refer the reader
to Chapter 6 of the monograph \cite{nualart06} for more details.

In the present article, we consider the parabolic Cauchy problem
\begin{eqnarray}
\label{parabolic-eq}
\frac{\partial u}{\partial t}(t,x)&=&\cL u(t,x)+\dot W(t,x), \quad t>0, x \in \bR^d \\
\nonumber
u(0,x) &=& 0, \quad x \in \bR^d,
\end{eqnarray}
and the hyperbolic Cauchy problem
 \begin{eqnarray}
\label{hyperbolic-eq}
\frac{\partial^2 u}{\partial t^2}(t,x)&=&\cL u(t,x)+\dot W(t,x), \quad t>0, x \in \bR^d \\
\nonumber
u(0,x) &=& 0,  \quad x \in \bR^d \\
\nonumber
 \frac{\partial u}{\partial t}(0,x) &=& 0,  \quad x \in
\bR^d,
\end{eqnarray}
where $\cL$ is a ``spatial operator'' (i.e. it acts only on the $x$ variable) given by the $L^2(\bR^d)$-generator of a $d$-dimensional
L\'evy process $X=(X_t)_{t \geq 0}$, and
$W$ is a Gaussian noise whose covariance is written formally as:
$$E[W(t,x)W(s,y)]=|t-s|^{2H-2}f(x-y),$$ for some index $H>1/2$
and some kernel $f$ (to be defined below).

The rigorous definition of the noise $\dot W$ is given in
Section \ref{prelim-section}. At this point, we should just mention that the covariance structure of the noise has two components:  a spatially-homogeneous component specified by the kernel $f$
(the example that we have in mind being the Riesz kernel
$f(x)=|x|^{-(d-\alpha)}$, with $0<\alpha<d$), and a temporal component
inherited from the fBm. This becomes clear once we realize that {\em if} $H>1/2$, $R_H(t,s)$ can be written as:
$$R_{H}(t,s)=\alpha_H\int_0^t \int_0^s |u-v|^{2H-2}du dv, \quad \mbox{with} \quad \alpha_H=H(2H-1).$$ (The case $H<1/2$ has to be treated differently and is not discussed here.)

The solution to problem (\ref{parabolic-eq}) (or
(\ref{hyperbolic-eq})) is understood in the mild-sense, and one of
the goals of the present article is to give a necessary and
sufficient condition for the existence of this solution, in terms
of the parameters $(H,f)$ of the noise, and the spatial operator
$\cL$. A similar problem has been considered in \cite{FKN09} and
\cite{FK10} in the case $H=1/2$ (which corresponds to the white
noise in time). This motivated us to examine the case $H>1/2$.

The case of the hyperbolic equation with spatial operator
$\cL=-(-\Delta)^{-\beta/2}, \linebreak \beta>0$, driven by a white noise in time was
examined in \cite{dalang-mueller03} and \cite{dalang-sanzsole05}.
In fact, these authors consider the much more difficult case of
the non-linear equation $\partial_{tt} u=\cL
u+\sigma(u)\dot{W}+b(u)$ with arbitrary initial conditions, and
Lipschitz continuous functions $\sigma$ and $b$. For the linear
equation, it turns out that the necessary and sufficient condition
for the existence of the solution is:
\begin{equation}
\label{Dalang-Mueller-cond}
\int_{\bR^d} \frac{1}{1+|\xi|^{\beta}}\mu(d\xi)<\infty,
\end{equation}
where the measure $\mu$ is the inverse Fourier transform of $f$ in $\cS'(\bR^d)$.

In the case of equations (\ref{parabolic-eq}) and
(\ref{hyperbolic-eq}) driven by a space-time white noise (i.e.
$H=1/2$ and $f=\delta_0$), the authors of \cite{FKN09} have shown
that the necessary and sufficient condition for the existence of
a random field solution is:
\begin{equation}
\label{FKN-cond}
\int_{\bR^d} \frac{1}{1+{\rm Re} \Psi(\xi)}\mu(d\xi)<\infty,
\end{equation}
where $\Psi(\xi)$ is the characteristic exponent of the underlying L\'evy process $X$. An important observation of \cite{FKN09} is that condition (\ref{FKN-cond}) can be extrapolated in a different context, being the condition which guarantees the existence of a local time $\int_0^t \delta_0(\bar{X}_s)ds$ of the symmetrization $\bar{X}$ of $X$.
This line of investigation was continued in \cite{FK10} in the case of the parabolic equation (\ref{parabolic-eq}) with white noise in time, but covariance kernel $f$ in space. Surprisingly, it is shown there that condition (\ref{FKN-cond}) is related not only to the existence of the ``occupation'' local time $L_t(f)=\int_{0}^{t}f(\bar{X}_s)ds$, but also to the potential theory of the process $\bar{X}$, when viewed as a Markov process.

In the present article, we carry out a similar program in the case
of the fractional noise in time. More precisely, after the
introduction of some background material in Section
\ref{prelim-section}, the article is split between the two
problems: Section \ref{parabolic-section} is dedicated to the
parabolic problem (\ref{parabolic-eq}), while Section
\ref{hyperbolic-section} treats the hyperbolic problem
(\ref{hyperbolic-eq}). For the parabolic problem, we discuss three
things: (i) the existence of a random field solution (Section 3.1); (ii) the maximal principle which gives the connection
with the potential theory of Markov processes (Section 3.2); (iii) the
relationship with the ``weighted'' intersection local time $L_{t,H}(f)$, defined by
$$L_{t,H}(f)=\alpha_H \int_0^t \int_0^t |r-s|^{2H-2} f(\bar{X}_r^1-\bar{X}_s^2)dr ds,$$
where $\bar{X}^1$ and $\bar{X}^2$ are two independent copies of
$\bar{X}$ (Section 3.3). For the hyperbolic problem, we only discuss the
existence of a random field solution in the case when
$\Psi(\xi)$ is real-valued (i.e. $X$ is symmetric).

Unlike the case of the white noise in time, it turns out that for the fractional noise, the conditions for the existence of the solution are different for the parabolic and hyperbolic problems. These conditions are:
\begin{equation}
\label{parabolic-cond} \int_{\bR^d} \left(\frac{1}{1+{\rm Re} \Psi(\xi)}\right)^{2H}\mu(d\xi)<\infty,
\end{equation}
in the parabolic case, respectively,
\begin{equation}
\label{cond-hyp-eq} \int_{\bR^d} \left(\frac{1}{1+{\rm Re}\Psi(\xi)}
\right)^{H+1/2}\mu(d\xi)<\infty,
\end{equation}
in the hyperbolic case. This phenomenon was
observed for the first time in \cite{BT10-SPA} for the wave and heat equations. 
As an application, we discuss
the case when $X$ is a $\beta$-stable process with $\beta \in (0,2]$, and hence $\Psi(\xi)=c_{\beta}|\xi|^{\beta}$. In this case, conditions
(\ref{parabolic-cond}) and (\ref{cond-hyp-eq}) turn out to be generalizations of (\ref{Dalang-Mueller-cond}).

The fact that the fractional noise induces a connection with the
weighted intersection local time was also used in
\cite{hu-nualart09}, in the case when $f=\delta_0$. In
\cite{BT09-JTP}, it is shown that the existence of the exponential
moment of $L_{t,H}(f)$ is closely related to the existence of the
(mild) solution of the heat equation with multiplicative noise.

We now introduce the notation used in the present article. The Fourier transform of a function $\varphi \in L^1(\bR^d)$ is defined by:
$$\cF \varphi (\xi)=\int_{\bR^d} e^{-i \xi \cdot x} \varphi(x) dx.$$
It is known that the Fourier transform can be extended to
$L^2(\bR^d)$ (see e.g. \cite{folland92}). Plancherel
theorem says that for any $\varphi, \psi \in L^2(\bR^d)$,
$$\int_{\bR^d} \varphi(x) \psi(x) dx= \frac{1}{(2\pi)^{d}} \int_{\bR^d} \cF \varphi(\xi) \overline{\cF \psi(\xi)}d\xi.$$

Let $\cS(\bR^d)$ be the Schwartz space of rapidly
decreasing infinitely differentiable functions on $\bR^d$.
A continuous linear functional on $\cS(\bR^d)$ is called a tempered
 distribution. Let $\cS'(\bR^d)$ be the space of
  tempered distributions. The Fourier transform $\cF S$
  of a functional $S \in \cS'(\bR^d)$ is defined by:
$$(\cF S, \varphi)=(S, \cF \varphi), \quad \forall \varphi \in \cS(\bR^d).$$

We refer the reader to \cite{schwartz66} for more details about the space $\cS'(\bR^d)$.

\section{Preliminaries}
\label{prelim-section}

In this section, we introduce some background material about the
Gaussian noise $W$ and the L\'evy process $X$.

\subsection{The Gaussian noise}

 As in \cite{dalang99}, we let $f:\bR^d \to [0,\infty]$ be a measurable locally integrable function (or a {\em kernel}). We assume that $f$ is the Fourier transform in $\cS'(\bR^d)$ of a tempered measure $\mu$, i.e.
\begin{equation}
\label{def-Fourier-mu}
\int_{\bR^d}f(x)\varphi(x)dx=\int_{\bR^d} \cF \varphi(\xi) \mu(d\xi), \quad \forall \varphi \in \cS(\bR^d).
\end{equation}

It follows that for any $\varphi,\psi \in \cS(\bR^d)$,
\begin{equation}
\label{conseq-def-mu}
\int_{\bR^d} \int_{\bR^d} \varphi(x)\psi(y)f(x-y)dx dy=\int_{\bR^d} \cF \varphi(\xi) \overline{\psi(\xi)}\mu(d\xi).
\end{equation}

Similarly to \cite{BT10-SPA}, we let $\cE$ be the set of
elementary functions of the form $$h(t,x)= \phi(t)\psi(x), \quad t
\geq 0, \ x \in \bR^d,$$ where $\phi$ is a linear combination of
indicator functions $1_{[0,a]}$ with $a>0$, and $\psi \in
\cS(\bR^d)$. We endow $\cE$ with the inner product:
$$\langle h_1, h_2 \rangle_{\cH \cP}=\alpha_H \int_{0}^{\infty}
\int_{0}^{\infty} \int_{\bR^{d}} \int_{\bR^d}|t-s|^{2H-2}
f(x-y)h_1(t,x)h_2(s,y)dx dy dr ds.$$ Let $\cH \cP$ be the
completion of $\cE$ with respect to the inner product $\langle
\cdot , \cdot \rangle_{\cH \cP}$. We note that the space $\cH \cP$
may contain distributions in both $t$ and $x$ variables.

We consider a zero-mean Gaussian process $\{W(h); h \in \cE\}$ with covariance
$$E(W(h_1)W(h_2))=\langle h_1, h_2 \rangle_{\cH \cP}.$$

The map $h \mapsto W(h)$ is an isometry between $\cE$ and the
Gaussian space of $W$, which can be extended to $\cH \cP$. This extension defines an
isonormal Gaussian process $W=\{W(h);h \in \cH \cP\}$.
We write $$
W(h)=\int_0^{\infty}\int_{\bR^d} h(t,x)W(dt,dx), \quad \mbox{for any} \ h \in \cH \cP.$$ This defines
the stochastic integral of an element $h \in \cH \cP$ with respect
to the noise $W$.

\subsection{The L\'evy process}

As in \cite{FK10}, we let $X=(X_t)_{t \geq 0}$ be a
$d$-dimensional L\'evy process with characteristic exponent
$\Psi(\xi)$. Hence, $X_0=0$ and for any $t>0$,
\begin{equation}
\label{Fourier-transform-Xt}
E(e^{-i \xi \cdot
X_t})=e^{-t \Psi(\xi)}, \quad \mbox{for all} \ \xi \in \bR^d.
\end{equation}

By the L\'evy-Khintchine formula (see e.g Theorem 8.1 of \cite{sato99}),
$$\Psi(\xi)=i\gamma \cdot \xi +\xi^T A \xi -\int_{\bR^d}(e^{-i \xi \cdot x}-1+i\xi \cdot x 1_{\{|x| \leq 1\}} )\nu(dx),$$
where $(\gamma,A,\nu)$ is the generating triplet of $X$. We
observe that for any $\xi \in \bR^d$,
\begin{eqnarray*}
{\rm Re} \Psi (\xi)&=& \xi^T A \xi +\int_{\bR^d}
[1-\cos (\xi \cdot x)] \nu(dx) \geq 0, \quad \mbox{and} \\
{\rm Im} \Psi (\xi)&=& \gamma \cdot \xi + \int_{\bR^d}[\sin (\xi
\cdot x)-\xi \cdot x 1_{\{|x| \leq 1\}}] \nu(dx).
\end{eqnarray*}


 $X$ is a homogenous Markov process with transition
probabilities:
$$Q_t(x;B)=P(X_{s+t} \in B|X_s=x)=P(X_t \in B-x),$$
for any $x \in \bR^d$ and Borel set $B \subset \bR^d$. Let
$(P_t)_{t \geq 0}$ be the associated semigroup, defined by:
$$(P_t \phi)
(x)=\int_{\bR^d}\phi(y)Q_t(x;dy)=E[\phi(x+X_t)],$$
for any bounded (or non-negative) measurable function $\phi: \bR^d \to \bR$.

Let $\cL$ be the $L^2(\bR^d)$-generator of $(P_t)_{t \geq 0}$,
defined by:
$$\cL \phi=\lim_{t \to 0}\frac{P_t \phi-\phi}{t} \quad
\mbox{in} \ L^2(\bR^d) \quad\mbox{(if it exists)}.$$ Note that
$\cL \phi$ exists if and only if $\phi \in {\rm Dom}(\cL)$, where
$${\rm Dom}(\cL)=\{\phi \in L^2(\bR^d);
(\cF \phi)\Psi \in L^2(\bR^d)\}$$ (see
p.16 of \cite{FK10}). Moreover, $\cL$ can be viewed as a
convolution operator with Fourier multiplier
 $\cF \cL=-\Psi$, i.e. for any $\phi \in {\rm Dom}(\cL)$ and $\xi \in \bR^d$,
$$\cF(\cL \phi)(\xi)=-\Psi(\xi)\cF \phi(\xi).$$

\section{The Parabolic Equation}
\label{parabolic-section}

In this section, we assume that the law of $X_t$ has a density denoted by $p_t$. This assumption allows us to identify the fundamental solution of $\partial_t u-\cL u=0$.

To see this, note that $Q_t(x; \cdot)$ has density $p_t( \cdot
-x)$, and $P_t \phi=\phi * \tilde p_t$, where $\tilde
p_t(x)=p_t(-x)$. Since the solution of the Kolmogorov's equation
$\partial_t u(t,x)=\cL u(t,x)$ with initial condition
$u(0,x)=u_0(x)$ is $$u(t,x)=(P_t u_0)(x)=\int_{\bR^d}
u_0(y)p_t(y-x)dy,$$
 it follows that the fundamental solution of $\partial_t u-\cL u=0$ is the function:
$$G(t,x)=p_t(-x), \quad t>0, x \in \bR^d.$$
From (\ref{Fourier-transform-Xt}), we obtain that:
\begin{equation}
\label{Fourier-G}
{\cF} G(t,\cdot)(\xi)=
\int_{\bR^d}e^{i \xi \cdot x} p_t(x)dx=E(e^{i \xi \cdot X_t})=e^{-t \overline{\Psi(\xi)}}.
\end{equation}

\subsection{Existence of the Random-Field Solution}

There are two equivalent ways of defining a random field solution for problem (\ref{parabolic-eq}).
Similarly to \cite{BT10-SPA}, one can say that the process $\{u(t,x); t \geq 0,x \in \bR^d\}$ defined by:
$$u(t,x)=\int_0^t \int_{\bR^d}G(t-s,x-y)W(ds,dy)$$
is a {\em random field solution} of (\ref{parabolic-eq}), provided that the stochastic integral above is well-defined, i.e. the integrand
\begin{equation}
\label{gtx-in-H}
\bR_{+} \times \bR^d \ni (s,y)\mapsto g_{tx}(s,y)=1_{[0,t]}(s)G(t-s,x-y) \quad \mbox{belongs to} \quad \cH \cP.
\end{equation}

 Since $g_{tx}$  satisfies conditions (i)-(iii) of Theorem 2.1 of \cite{BT10-SPA}, to check that $g_{tx} \in \cH \cP$, it suffices to prove that:
$$I_t:=\alpha_H \int_{\bR^d} \int_{0}^{\infty} \int_{0}^{\infty} \cF g_{tx}(r,\cdot)(\xi) \overline{\cF g_{tx}(s,\cdot)(\xi)}|r-s|^{2H-2} dr ds \mu(d\xi)<\infty.$$
In this case,
$$E|u(t,x)|^2=E|W(g_{tx})|^2=
\|g_{tx}\|_{\cH \cP}^2=I_t.$$

Note that by (\ref{Fourier-G}), $I_t=\int_{\bR^d}N_t(\xi)\mu(d
\xi)$, where
$$N_t(\xi)=\alpha_H \int_0^t \int_0^t e^{-r \Psi(\xi)} e^{-s \overline{\Psi(\xi)}}|r-s|^{2H-2}dr ds.$$
Therefore, the question about the existence of a
random field solution of (\ref{parabolic-eq}) reduces to finding suitable upper and
 lower bounds for $N_t(\xi)$. 

Alternatively, the authors of \cite{FKN09} suggest a different
method for defining the random field solution of
(\ref{parabolic-eq}), which has the advantage that can be applied
also to the hyperbolic problem (\ref{hyperbolic-eq}) (for which
one cannot identify the fundamental solution $G$). Since this is
the method that we use in the present article, we explain it
below.

We say that $\{u(t,\varphi);t \geq 0, \varphi \in \cS(\bR^d)\}$ is a {\em the weak solution} of (\ref{parabolic-eq}) if
$$u(t,\varphi)=\int_0^t \int_{\bR^d}
\left(\int_{\bR^d}G(t-s,x-y)\varphi(x) dx \right) W(ds,dy).$$ Note
that the stochastic integral above is well-defined (as a random
variable in $L^2(\Omega)$) if and only if the integrand $$\bR_{+}
\times \bR^d \ni (s,y) \mapsto
h_{t,\varphi}(s,y)=1_{[0,t]}(s)(\varphi * p_{t-s})(y) \quad
\mbox{belongs to $\cH \cP$}.$$

Since $h_{t,\varphi}$ satisfies
the conditions (i)-(iii) of Theorem 2.1 of \cite{BT10-SPA}, to check that $h_{t,\varphi} \in
\cH \cP$ it suffices to show that
$$I_{t,\varphi}:=\alpha_H \int_{\bR^d} \int_0^{\infty} \int_0^{\infty}|r-s|^{2H-2} \cF h_{t,\varphi}(r,\cdot)(\xi)
\overline{\cF h_{t,\varphi}(s,\cdot)(\xi)}dr ds \mu(d\xi)<\infty.$$
In this case, we have
$$E|u(t,\varphi)|^2=E|W(h_{t,\varphi})|^2=
\|h_{t,\varphi}\|_{\cH \cP}^2=I_{t,\varphi}.$$

Since both $\varphi$ and $p_{t-s}$ are in $L^1(\bR^d)$,
$$\cF(\varphi * p_{t-s})(\xi)=\cF \varphi(\xi) \cF p_{t-s}(\xi)=\cF
\varphi(\xi) e^{-(t-s) \Psi(\xi)},$$ and hence
$$I_{t,\varphi}=\int_{\bR^d}N_t(\xi)|\cF \varphi(\xi)|^2 \mu(d\xi).$$

Using the trivial bound $N_t(\xi) \leq t^{2H}$ (which is obtained using the fact that $|e^{-s \Psi(\xi)}|=e^{-s {\rm Re} \Psi(\xi)} \leq 1$ for all $s>0$ and $\xi \in \bR^d$), we get: $$I_{t,\varphi}\leq t^{2H} \int_{\bR^d}|\cF \varphi(\xi)|^2 \mu(d\xi)<\infty \quad \mbox{for all} \ \varphi \in \cS(\bR^d).$$
This shows that $u(t,\varphi)$ is a well-defined random variable in $L^2(\Omega)$.

We continue to explain the method of \cite{FKN09}. We endow $\cS(\bR^d)$ with the inner product:
$$\langle \varphi,\psi \rangle_t=E(u(t,\varphi)u(t,\psi)).$$

We denote by $\| \cdot \|_t$ the norm induced by the inner product $\langle \cdot, \cdot \rangle_t$, i.e.
\begin{equation}
\label{def-norm-2bars}
\|\varphi \|_t^2=E|u(t,\varphi)|^2 =\int_{\bR^d} N_t(\xi)|\cF \varphi(\xi)|^2 \mu(d\xi).
\end{equation}

Let $M_t$ be the completion of $\cS(\bR^d)$ with respect to 
$\langle \cdot, \cdot \rangle_t$ and $M=\cap_{t>0}M_t$. We say that (\ref{parabolic-eq}) has a {\em random field solution} if and only if
\begin{equation}
\label{delta-x-in-M}
\delta_x \in M \quad \mbox{for all} \ x \in \bR^d.
\end{equation}
The random field solution is defined by $\{u(t,x)=u(t,\delta_x); t \geq 0,x \in \bR^d\}$.

To prove (\ref{delta-x-in-M}), we introduce the space $\cZ=\cap_{t>0}\cZ_t$, where $\cZ_t$ is the
completion of $\cS(\bR^d)$ with respect to the inner product $[\cdot,\cdot]_t$ defined by:
$$[\varphi, \psi]_t= \int_{\bR^d} \left(\frac{1}{1/t+ {\rm Re} \Psi(\xi)} \right)^{2H} \cF \varphi(\xi) \overline{\cF \psi(\xi)}\mu(d\xi)=:\cE(t; \varphi,\psi)$$

We denote by $|\| \cdot \||_t$ the norm induced by the inner product $[\cdot, \cdot]_t$, i.e.
\begin{equation}
\label{def-norm-3bars}
|\| \varphi \||_t  =\int_{\bR^d} \left(\frac{1}{1/t +  {\rm Re} \Psi(\xi)} \right)^{2H} |\cF \varphi(\xi)|^2\mu(d\xi)=:\cE(t; \varphi).
\end{equation}

By Lemma \ref{ineq-cE} (Appendix A), for any $s,t>0$, there exist some positive constants $c_1(s,t)$ and $c_2(s,t)$ such that
for any $\varphi \in \cS(\bR^d)$,
$$c_1(s,t)^{2H} \cE(s;\varphi) \leq \cE(t; \varphi) \leq c_2(s,t)^{2H} \cE(s; \varphi).$$
Therefore, the norms $|\| \cdot \||_t$ and $|\| \cdot \||_s$ are equivalent and $\cZ_t =\cZ_s=\cZ$.

The idea for proving (\ref{delta-x-in-M}) is to show that any norm $\| \cdot \|_t$ is equivalent to a norm $|\| \cdot \||_{\rho(t)}$, for a certain bijective function $\rho:\bR_{+} \to \bR_{+}$. From this, one infers that $M_t=\cZ_{\rho(t)}$ for any $t>0$, and hence $M=\cZ=\cZ_1$. Condition (\ref{delta-x-in-M}) becomes $\delta_x \in \cZ_1$ for all $x \in \bR^d$, for which one can find a natural necessary and sufficient condition (see Theorem \ref{nec-suf-delta-in-Z1} below).
In the case of the parabolic problem (\ref{parabolic-eq}), it turns out that
$\rho(t)=t$. (We will see in Section \ref{hyperbolic-section} that for the hyperbolic problem (\ref{hyperbolic-eq}), $\rho(t)=t^2$.)

The next theorem is the main result of the present section, and
gives the desired upper and lower bounds for $N_t(\xi)$. Unfortunately,
for the lower bound, we had to introduce an additional condition of boundedness on the
ratio between the imaginary part and the real part of the characteristic exponent $\Psi(\xi)$.
A similar difficulty has been encountered in \cite{khoshnevisan-xiao09} for
obtaining a lower bound for the ``sojourn operator''. Our condition (\ref{cond-Im-Re})
is similar to condition (3.3) of \cite{khoshnevisan-xiao09}, and is trivially
 satisfied when $\Psi$ is real-valued.

We use the following inequality: there exists a constant $b_H>0$, such that
\begin{equation}
\label{MMV-ineq} \alpha_H \int_{\bR}
\int_{\bR} |\varphi(r)| |\varphi(s)||r-s|^{2H-2}dr ds \leq b_H^2 \left(
\int_{\bR} |\varphi(s)|^{1/H}ds \right)^{2H}
\end{equation}
for any 
$\varphi \in L^{1/H}(\bR)$. This inequality  was proved in
\cite{MMV01} and is a consequence of the Littlewood-Hardy
inequality.

For complex-valued functions $\varphi$, we define:
$$ \|\varphi \|_{\cH(0,t)}^2:=\alpha_H \int_0^t
\int_0^t \varphi(r) \overline{\varphi(s)}|r-s|^{2H-2}dr ds =\|{\rm Re} \varphi \|_{\cH(0,t)}^2+ \|{\rm Im} \varphi\|_{\cH(0,t)}^{2}.$$

\begin{theorem}
\label{bounds-Nt-xi} 
For any $t>0$ and $\xi \in \bR^d$,
$$N_t(\xi) \leq  C_H \left(\frac{1}{1/t+ {\rm Re} \Psi(\xi)} \right)^{2H},$$
where $C_H =H^{2H} b_H^2 e^{2}$.
If in addition, there exists a constant $K>0$ such that:
\begin{equation}
\label{cond-Im-Re}
|{\rm Im} \Psi (\xi)| \leq K {\rm Re} \Psi(\xi), \quad \forall \xi \in \bR^d.
\end{equation}
then, for any $t>0$ and $\xi \in \bR^d$,
$$N_t(\xi) \geq C_{H,K} \left(\frac{1}{1/t+ {\rm Re} \Psi(\xi)} \right)^{2H},$$
where $C_{H,K}$ is a positive constant depending on $H$ and $K$.
\end{theorem}

\noindent {\bf Proof:} For the upper bound, we note that $N_t(\xi)$ can be written as
$$N_t(\xi)=\alpha_H \int_0^{t} \int_0^t e^{-r {\rm Re} \Psi(\xi)} e^{-s {\rm Re} \Psi(\xi)} |r-s|^{2H-2}\cos[(r-s){\rm Im} \Psi(\xi)] dr ds.$$

Using the fact that $|\cos x| \leq 1$ and $e^{-r {\rm Re} \Psi(\xi)} \leq e^{t/\lambda} e^{-r(1/\lambda+{\rm Re} \Psi(\xi))}$ for any $r \in [0,t]$, we get:
$$N_t(\xi) \leq e^{2t/\lambda} \alpha_H \int_0^{t} \int_0^t e^{-r(1/\lambda +{\rm Re} \Psi(\xi))} e^{-s(1/\lambda +{\rm Re} \Psi(\xi))} |r-s|^{2H-2}dr ds $$

By (\ref{MMV-ineq}), it follows that:
$$N_t(\xi) \leq  e^{2t/\lambda} b_H^2 \left(\int_0^t e^{-r(1/\lambda+{\rm Re} \Psi(\xi))/H}dr \right)^{2H} \leq b_H^2  e^{2t/\lambda}\left(\frac{H}{1/\lambda+ {\rm Re} \Psi(\xi)}\right)^{2H}.$$
The conclusion follows by taking $\lambda=t$.

For the lower bound, suppose first that $t {\rm Re} \Psi(\xi) \leq a$, for some constant $a=a_K \in (0,1)$
such that $Ka < \pi/2$.
By (\ref{cond-Im-Re}), $$t|{\rm Im} \Psi(\xi)| \leq K t {\rm Re} \Psi (\xi) \leq Ka <\frac{\pi}{2}.$$

Using the fact that $e^{-x} \geq 1-x$  for $x>0$, we obtain: for any $r \in [0,t]$,
$$e^{-r {\rm Re} \Psi(\xi)} \geq 1-r {\rm Re} \Psi (\xi)\geq 1-a.$$
 Since $\cos$ is decreasing on the interval $[0,\frac{\pi}{2}]$, for any $0<s<r<t$, $$\cos[(r-s)|{\rm Im} \Psi(\xi)|] \geq \cos [t|{\rm Im} \Psi(\xi)|] \geq \cos (Ka)>0.$$
Therefore,
\begin{eqnarray*}
N_t(\xi)&=&2\alpha_H\int_0^{t} \int_0^r e^{-r {\rm Re} \Psi(\xi)} e^{-s {\rm Re} \Psi(\xi)} (r-s)^{2H-2}\cos[(r-s)|{\rm Im} \Psi(\xi)|] ds dr \\
& \geq &  (1-a)^2 \cos (Ka) 2\alpha_H \int_0^{t} \int_0^r (r-s)^{2H-2} ds dr \\
&=&  (1-a)^2 \cos (Ka) t^{2H} \geq (1-a)^2 \cos (Ka) \left(\frac{1}{1/t+{\rm Re} \Psi(\xi)} \right)^{2H},
\end{eqnarray*}
where for the last inequality we used the fact that $t \geq \frac{1}{t^{-1}+{\rm Re} \Psi(\xi)}$.

Suppose next that $t {\rm Re} \Psi(\xi) \geq a$. Note that
$$N_t(\xi)=\|e^{- \cdot  \Psi(\xi)} \|_{\cH(0,t)}^2.$$

Using Lemma \ref{H-norm-exp} (Appendix B) for expressing the $\cH(0,t)$-norm of the exponential function
in the spectral domain, we obtain:
$$N_t(\xi)=c_{H} \int_{\bR} \frac{\sin^2[(\tau+ {\rm Im} \Psi(\xi))t]+\{e^{-t {\rm Re} \Psi(\xi)}-\cos[(\tau+{\rm Im} \Psi(\xi))t]\}^2}{[{\rm Re} \Psi(\xi)]^2 + [\tau+ {\rm Im} \Psi(\xi)]^2} |\tau|^{-(2H-1)} d\tau.$$

We denote
$$T=t {\rm Re} \Psi(\xi) \quad \mbox{and} \quad b=\frac{{\rm Im} \Psi(\xi)}{{\rm Re} \Psi(\xi)}.$$
Using the change of variable $\tau'=\tau/{\rm Re} \Psi(\xi)$, we obtain that:
\begin{equation}
\label{new-form-Ntxi-par}
N_t(\xi)=\frac{c_H}{[{\rm Re} \Psi(\xi)]^{2H}} \int_{\bR}\frac{|\tau|^{-(2H-1)}}{1+(\tau+b)^2}[f_T^2 (\tau)+g_T^{2}(\tau)] d\tau,
\end{equation}
where
$f_T(\tau)=\sin [(\tau+b)T]$ and $g_T(\tau)=e^{-T}-\cos[(\tau+b)T]$.

From the proof of Lemma \ref{H-norm-exp} (Appendix B), we know that:
$$\frac{1}{1+(\tau+b)^2}[f_T^2 (\tau)+g_T^{2}(\tau)]=|\cF_{0,T} \varphi(\tau)|^2,$$ where $\varphi(x)=e^{-x(1+ib)}$.

Let $\rho>K$ be positive constant whose value will be specified later. Since the integrand of (\ref{new-form-Ntxi-par}) is non-negative, the integral can be bounded below by the integral over the region $|\tau| \leq \rho$. In this region, $|\tau|^{-(2H-1)} \geq \rho^{-(2H-1)}$. We obtain:
\begin{equation}
\label{LB-Nt-heat-step1}
N_t(\xi) \geq  \frac{c_H \rho^{-(2H-1)}}{[{\rm Re} \Psi(\xi)]^{2H}}
\left(I(T) - \int_{|\tau| \geq \rho}\frac{1}{1+(\tau+b)^2}[f_T^2 (\tau)+g_T^{2}(\tau)]  d\tau \right),
\end{equation}
where
\begin{equation}
\label{Planch-id-heat}
I(T):=\int_{\bR}\frac{1}{1+(\tau+b)^2}[f_T^2 (\tau)+g_T^{2}(\tau)] d\tau=2 \pi \int_0^T |e^{-x(1+ib)}|^2 dx=\pi(1-e^{-2T}),
\end{equation}
by Plancherel's theorem.

Using (\ref{Planch-id-heat}), we obtain the lower bound:
\begin{equation}
\label{LB-Nt-heat-step2}
I(T) \geq \pi(1-e^{-2a}), \quad \mbox{since} \quad T \geq a.
\end{equation}

To find an upper bound for the second integral on the right-hand side of (\ref{LB-Nt-heat-step1}), we use the fact that:
$$f_T^2(\tau)+g_T^2(\tau) \leq 5, \quad \forall \tau \in \bR.$$

It follows that:
\begin{equation}
\label{LB-Nt-heat-step3}
\int_{|\tau| \geq \rho} \frac{f_T^{2}(\tau)+g_T^{2}(\tau)}{1+(\tau+b)^2}d\tau \leq \int_{|\tau| \geq \rho}\frac{5}{(\tau+b)^2}d\tau =
\frac{10\rho}{\rho^2-b^2} \leq \frac{10\rho}{\rho^2-K^2},
\end{equation}
since $|b| \leq K$ (by (\ref{cond-Im-Re})). We choose $\rho=\rho_K$ large enough such that
$$C_K:=\pi(1-e^{-2a})-\frac{10\rho}{\rho^2-K^2}>0.$$

Using (\ref{LB-Nt-heat-step1}), (\ref{LB-Nt-heat-step2}) and  (\ref{LB-Nt-heat-step3}), we obtain:
$$N_t(\xi) \geq  C_K \frac{c_H \rho^{-(2H-1)}}{[{\rm Re} \Psi(\xi)]^{2H}}  \geq C_{K}  c_H \rho^{-(2H-1)}  \left(\frac{1}{1/t+{\rm Re \Psi(\xi)}} \right)^{2H}.$$

The conclusion follows, letting
$$C_{H,K}=\min\left\{(1-a)^2 \cos (Ka), C_{K}  c_H \rho^{-(2H-1)} \right\}.$$

$\Box$

The following result is an immediate consequence of Theorem \ref{bounds-Nt-xi}.

\begin{corollary}
\label{corol-bounds-Nt-xi}
a) For any $t>0, \varphi \in \cS(\bR^d)$,
$$E|u(t,\varphi)|^2 \leq C_H \cE(t;\varphi),$$
where $C_H=H^{2H} b_H^2 e^{2}$. Hence, $M_t \supset \cZ_t$ for all
$t>0$, and $M \supset \cZ$.

b) If (\ref{cond-Im-Re}) holds, then
for any $t>0, \varphi \in \cS(\bR^d)$,
$$E|u(t,\varphi)|^2 \geq C_{H,K} \cE(t;\varphi),$$
where $C_{H,K}$ is a positive constant depending on $H$ and $K$. Hence, $M_t=\cZ_t$ for all $t>0$, and $M=\cZ$.
\end{corollary}

\noindent {\bf Proof:} We use Theorem \ref{bounds-Nt-xi} and the definitions (\ref{def-norm-2bars}) and (\ref{def-norm-3bars}) of the norms $\|\cdot \|_t$, respectively $|\| \cdot \||_t$. $\Box$

\vspace{3mm}

The next result gives the necessary and sufficient condition for $\delta_x \in \cZ_1$ for all $x \in \bR^d$. 

\begin{theorem}
\label{nec-suf-delta-in-Z1}
In order that $\delta_x \in \cZ_1$ for all $x \in \bR^d$, it is necessary and sufficient that
condition (\ref{parabolic-cond}) holds.
\end{theorem}

\noindent {\bf Proof:} Suppose first that (\ref{parabolic-cond}) holds. To show that $\delta_x \in \cZ_1$ for all $x \in \bR^d$, we use an argument similar to the proof of Theorem 2 of \cite{dalang99}.

Let $\cZ_0$ be the set Schwartz distributions $\varphi$ such that $\cF \varphi$ is a function and
$$|\|\varphi \||_1:=\int_{\bR^d} \left(\frac{1}{1+{\rm Re} \Psi(\xi)} \right)^{2H}|\cF \varphi(\xi)|^2 \mu(d\xi)<\infty.$$

Note that $\cS(\bR^d) \subset \cZ_0$ and the definition of $|\| \cdot \| |_1$ agrees on $\cS(\bR^d)$ with the one given by (\ref{def-norm-3bars}). Therefore, to show that a distribution $\varphi \in \cZ_0$ is in $\cZ_1$, it suffices to show that there exists a sequence $(\varphi_n)_{n \geq 1} \subset \cS(\bR^d)$ such that $|\|\varphi_n-\varphi \||_1 \to 0$. We apply this to $\varphi=\delta_x$. In this case, $\cF \varphi(\xi)=e^{-i \xi \cdot x}$, $|\cF \varphi (\xi)|=1$ for all $\xi \in \bR^d$, and
$|\| \varphi \||_1$ coincides with the integral of (\ref{parabolic-cond}).

Let $\varphi_n=\varphi * \phi_n \in \cS(\bR^d)$, where $\phi_n(x)=n^d \phi(nx)$ and $\phi \in \cS(\bR^d)$ is such that $\phi \geq 0$ and $\int_{\bR^d}\phi(x)dx=1$.
Then $\cF \varphi_n(\xi)=\cF \varphi (\xi) \cF \phi_n(\xi)$ and
\begin{eqnarray*}
|\|\varphi_n-\varphi \||_1 &=& \int_{\bR^d}\left(\frac{1}{1+ {\rm Re} \Psi(\xi)} \right)^{2H} |\cF \varphi_n (\xi) -\cF \varphi (\xi)|^2 \mu(d\xi) \\
&=& \int_{\bR^d}\left(\frac{1}{1+ {\rm Re} \Psi(\xi)} \right)^{2H} |\cF \phi_n (\xi) -1|^2 \mu(d\xi) \to 0,
\end{eqnarray*}
by the Dominated Convergence Theorem, since $|\cF \phi_n(\xi)| \leq 1$ for all $\xi \in \bR^d$.

For the reverse implication, suppose that $\delta_x \in \cZ_1$ for all $x \in \bR^d$. To show that (\ref{parabolic-cond}) holds, one can use the same argument as in the proof of Lemma 4.2 of \cite{FK10}. We omit the details. $\Box$

\vspace{3mm}

The following result concludes our discussion about the existence of a random-field solution.

\begin{theorem}
\label{th-existence-parabolic}
{\rm (Existence of Solution in the Parabolic Case)}

a) If (\ref{parabolic-cond}) holds,
then equation (\ref{parabolic-eq}) has a random-field solution.

b) Suppose that (\ref{cond-Im-Re}) holds.  If (\ref{parabolic-eq}) has a random-field solution, then (\ref{parabolic-cond}) holds.
\end{theorem}

\noindent {\bf Proof:} a) Suppose that (\ref{parabolic-cond}) holds. By Theorem \ref{nec-suf-delta-in-Z1}, $\delta_x \in \cZ_1$ for all $x \in \bR^d$.  By Corollary \ref{corol-bounds-Nt-xi}.a), $\cZ_1=\cZ \subset M$. Hence (\ref{delta-x-in-M}) holds.

b) Suppose that (\ref{delta-x-in-M}) holds. By Corollary \ref{corol-bounds-Nt-xi}.b), $M=\cZ=\cZ_1$. Hence $\delta_x \in \cZ_1$ for all $x \in \bR^d$.  By Theorem \ref{nec-suf-delta-in-Z1}, (\ref{parabolic-cond}) holds. $\Box$

\begin{example}
\label{stable-ex}{\rm (Stable processes)}
{\rm Suppose that $\cL=-(-\Delta)^{\beta/2}$ for $\beta \in (0,2]$. Then $(X_t)_{t \geq 0}$ is a rotation invariant strictly $\beta$-stable process on $\bR^d$, and $\Psi(\xi)=c_{\beta}|\xi|^{\beta}$ (see Theorem 14.14  in \cite{sato99} and Example 30.6 in \cite{sato99}).
It can be shown that $(X_t)_{t \geq 0}$ is subordinate to the Brownian motion on $\bR^d$ by a strictly ($\beta/2$)-stable subordinator (see Example 32.7 of \cite{sato99}).

In this case, condition (\ref{parabolic-cond}) becomes:
\begin{equation}
\label{cond-stable-proc}
\int_{\bR^d} \left(\frac{1}{1+|\xi|^{\beta}} \right)^{2H}\mu(d\xi)<\infty.
\end{equation}

We consider two kernels:

(i) $f(x)=c_{\alpha,d}|x|^{-(d-\alpha)}$ for $0<\alpha<d$. In this case, $\mu(d\xi)=|\xi|^{-\alpha}d\xi$ (see p.117 of \cite{stein70}), and condition (\ref{cond-stable-proc}) is equivalent to
\begin{equation}
\label{cond-H-alpha-beta}
2H \beta>d-\alpha.
\end{equation}

(ii) $f(x)=\prod_{i=1}^{d}(\alpha_{H_i}|x_i|^{2H_i-2})$. In this case, $\mu(d\xi)=\prod_{i=1}^{n}c_{H_i}|\xi_i|^{-(2H_i-1)}$, and condition (\ref{cond-stable-proc}) is equivalent to
$$2H \beta >d-\sum_{i=1}^{d}(2H_i-1).$$
}
\end{example}

\begin{remark} {\rm (Fractional Powers of the Laplacian)}
\label{stable-rem}
{\rm As in \cite{dalang-mueller03} and \cite{dalang-sanzsole05}, we can consider also the case $\cL=-(-\Delta)^{\beta/2}$ for arbitrary $\beta>0$, even if there is no corresponding L\'evy process whose generator is $\cL$. Note that the fundamental solution $G$ of $\partial_t u- \cL u=0$ exists and satisfies:
$$\cF G (t,\xi)=\exp(-c_{\beta}t|\xi|^{\beta}).$$
Using Theorem 2.1 of \cite{BT10-SPA} and estimates similar to those given by Theorem \ref{bounds-Nt-xi} above, one can show that a random field solution of (\ref{parabolic-eq}) (in the sense of \cite{BT10-SPA}) exists if and only if (\ref{cond-stable-proc}) holds.
}
\end{remark}

\subsection{A Maximum Principle}

Throughout this section, we assume that $p_t \in L^2(\bR^d)$ for all $t>0$, and
\begin{equation}
\label{mu-has-density}
\mbox{$\mu$ has a (non-negative) density $g$,}
\end{equation}
i.e.  $f$ is a kernel of positive type (see Definition 5.1 of \cite{khoshnevisan-xiao09}).

We consider the symmetric L\'evy process $\bar{X}=(\bar{X}_t)_{t \geq 0}$ defined by:
$$\bar{X}_t:=X_t-\tilde{X}_t,$$
where $(\tilde{X}_t)_{t \geq 0}$ is an independent copy of $(X_t)_{t \geq 0}$.
We denote by $(\bar{P}_t)_{t \geq 0}$ the semigroup of $(\bar{X}_t)_{t \geq 0}$, i.e.
$$(\bar{P}_t \phi)(x)=\int_{\bR^d}\phi(y)\bar{p}_t(x-y)dy,$$
where $\bar{p}_t=p_t*\tilde{p}_t$. From (\ref{Fourier-G}), it follows that $\cF \bar{p}_{t}(\xi)=e^{-2t {\rm Re} \Psi(\xi)}$.

Let $(\bar{R}_{\alpha})_{\alpha>0}$ be the resolvent of $(\bar{P}_t)_{t \geq 0}$, i.e.
$$(\bar{R}_{\alpha}\phi)(x)=\int_{0}^{\infty}e^{-\alpha s} (\bar{P}_s \phi)(x)ds.$$


The following maximum principle has been obtained recently in \cite{FK10}:  
\begin{equation}
\label{max-principle-FK}
(\bar{R}_{\alpha}f)(0)=\sup_{x \in \bR^d}(\bar{R}_{\alpha}f)(x)=\Upsilon(\alpha):=\int_{\bR^d} \frac{1}{\alpha+2{\rm Re} \Psi(\xi)}\mu(d\xi).
\end{equation}

\begin{remark}
\label{remark-kernel-of-positive-type}
{\rm Recall that $f=\cF g$ in $\cS'(\bR^d)$.
The authors of \cite{FK10} work with the Fourier transform $\cF f$ instead of $g$, which introduces an additional factor $(2\pi)^{-d}$. To see this, note that by the Fourier inversion theorem on $\cS(\bR^d)$, 
relation (\ref{def-Fourier-mu}) becomes: for any $\phi \in \cS(\bR^d)$,
$$\int_{\bR^d}  \phi(\xi) g(\xi) d\xi=\frac{1}{(2\pi)^d} \int_{\bR^d}f(x)\cF \phi(x)dx.$$
This shows that $g=(2\pi)^{-d} \cF f$ in $\cS'(\bR^d)$.
}
\end{remark}

Note that $\Upsilon(\alpha)<\infty$ for all $\alpha>0$ if and only if $\Upsilon(\alpha)<\infty$ for some $\alpha>0$. An important consequence of (\ref{max-principle-FK}) (combined with the results of \cite{dalang99}) is that the potential-theoretic condition:
$$(\bar{R}_{\alpha}f)(0)<\infty \quad \mbox{for all} \quad \alpha>0$$
is necessary and sufficient for the existence of a random field solution of (\ref{parabolic-eq}), when the Gaussian noise $W$ is white in time (i.e. $H=1/2$). 

\vspace{3mm}

In the present article, we develop a maximal principle similar to (\ref{max-principle-FK}), which has a connection with the existence of a random field solution of (\ref{parabolic-eq}), when the noise $W$ is fractional in time.

We define the following ``fractional analogue'' of the
resolvent operator:
$$(\bar{R}_{\alpha,H} \phi)(x)=\alpha_H \int_0^{\infty} \int_0^{\infty}|r-s|^{2H-2}e^{-\alpha (r+s)} (\bar{P}_{r+s}\phi)(x)dr ds,$$
and we let
$$\Upsilon_{H}(\alpha):=\alpha_H \int_{\bR^d} \int_0^{\infty} \int_0^{\infty} |r-s|^{2H-2} e^{-(\alpha+ 2{\rm Re}\Psi(\xi)) (r+s)}  dr ds  \mu(d\xi).$$

As in \cite{khoshnevisan-xiao09}, we assume that $f$ satisfies the following condition:
\begin{equation}
\label{cond-on-f}
\mbox{$f(x)<\infty$ if and only if $x \not =0$.}
\end{equation}
Under this condition, the following harmonic-analysis result holds:
\begin{equation}
\label{harmonic-result}
\int_{\bR^d} \int_{\bR^d} \varphi(x)\psi(y)f(x-y)dx dy=\int_{\bR^d} \cF \varphi(\xi)\overline{\cF \psi(\xi)}g(\xi)d\xi,
\end{equation}
for any non-negative functions $\varphi,\psi \in L^1(\bR^d)$ (see Lemma 5.6 of \cite{khoshnevisan-xiao09}).

\begin{theorem}{\rm (A maximum principle)}
\label{max-principle} If (\ref{cond-on-f}) holds, then for any $\alpha>0$,
$$(\bar{R}_{\alpha,H}f)(0)=\sup_{x \in \bR^d} (\bar{R}_{\alpha,H}f)(x)=\Upsilon_{H}(\alpha).$$
\end{theorem}

The proof of Theorem \ref{max-principle} follows from Lemma \ref{max-principle-step3} and Lemma \ref{max-principle-step4} below. Before this, we need some intermediate results.
Let $C_0(\bR^d)$ be the space of continuous functions which vanish at infinity.

\begin{lemma}
\label{max-principle-step1}
For any $\phi \in \cS(\bR^d)$, we have: \\
a) $f*\phi \in C_0(\bR^d)$ and $\bar{R}_{\alpha,H}(f*\phi) \in C_0(\bR^d)$ for any $\alpha>0$;\\
b)  $f*\phi \in L^2(\bR^d)$ and $\cF (f*\phi)(\xi)=(2\pi)^d \cF \phi(\xi) g(\xi)$.
\end{lemma}

\noindent {\bf Proof:} a) Since $f$ is tempered (i.e. $f(x) \leq C (1+|x|)^k$ for all $x \in \bR^d$, for some $k \geq 0,C>0$),
the function $f*\phi$ is well-defined. By (\ref{def-Fourier-mu}) and (\ref{mu-has-density}),
\begin{equation}
\label{formula-f*phi}
(f*\phi)(x)=\int_{\bR^d}f(y)\phi(x-y)dy=\int_{\bR^d}e^{-i \xi \cdot x} \overline{\cF \phi(\xi)}g(\xi)d\xi.
\end{equation}
Since $g$ is tempered, $(\overline{\cF \phi}) g \in L^p(\bR^d)$ for any $p \geq 1$. By Riemann-Lebesgue lemma, $f*\phi \in C_0(\bR^d)$. Finally, we note that $\bar{R}_{\alpha,H} :C_0(\bR^d) \to C_0(\bR^d)$.

b) By (\ref{formula-f*phi}), $f*\phi=\cF h$, where $h:=(\overline{\cF \phi}) g \in L^2(\bR^d)$. Hence, $f*\phi \in L^2(\bR^d)$. By the Fourier inversion formula in $L^2(\bR^d)$ (see e.g. p.222 of \cite{folland92}), $$\cF (f*\phi)(\xi)=(2\pi)^{d}\overline{h(\xi)}=(2\pi)^d  {\cF \phi}(\xi)g(\xi).$$
$\Box$

\begin{lemma}
\label{max-principle-step2}
For any $\phi \in \cS(\bR^d)$ and $x \in \bR^d$,
$$(\bar{R}_{\alpha,H}(f*\phi))(x)=\alpha_H \int_{\bR^d}  e^{-i\xi \cdot x}\overline{\cF \phi(\xi)} \int_{\bR_{+}^2} |r-s|^{2H-2} e^{-(\alpha+2{\rm Re} \Psi(\xi))(r+s)}dr ds \mu(d\xi).$$
Consequently, for any $\phi \in \cS(\bR^d)$ with $\|\phi\|_1=1$,
$$|(\bar{R}_{\alpha,H}(f*\phi))(x)| \leq \Upsilon_{H}(\alpha), \quad \forall  x \in \bR^d.$$
\end{lemma}

\noindent {\bf Proof:} Since $\bar{P}_{r+s}=\bar{P}_r \bar{P}_s$, we have:
\begin{equation}
\label{calcul-P-st}
(\bar{P}_{r+s}(f*\phi))(x)=\int_{\bR^d}\int_{\bR^d} (f*\phi)(y-z)\bar{p}_r(x-y) \bar{p}_s(z) dy dz.
\end{equation}
Using Lemma \ref{version-Plancherel} (Appendix C)
with $\varphi=f*\phi$, $\psi_1=\bar{p}_r(x-\cdot)$ and $\psi_2=\bar{p}_s$,
$$(\bar{P}_{r+s}(f*\phi))(x)=\frac{1}{(2\pi)^d} \int_{\bR^d} e^{-i \xi \cdot x} \overline{\cF \bar{p}_r(\xi)} \ \overline{\cF \bar{p}_s(\xi)} \ \overline{\cF(f*\phi)(\xi)}d\xi.$$
The result follows using Lemma \ref{max-principle-step1}.b), the fact that $\cF \bar{p}_{r}(\xi)=e^{-2r {\rm Re} \Psi(\xi)}$, and Fubini's theorem. $\Box$

\begin{lemma}
\label{max-principle-step3} For any $\alpha>0$,
$$\Upsilon_{H}(\alpha)=\sup_{x \in \bR^d}(\bar{R}_{\alpha,H}f)(x)=\limsup_{x \to 0} (\bar{R}_{\alpha,H}f)(x).$$
\end{lemma}

\noindent {\bf Proof:} The proof is similar to Proposition 3.5 of \cite{FK10}.
By Fatou's lemma and Lemma \ref{max-principle-step2}, for any $x \in \bR^d$,
$$(\bar{R}_{\alpha,H}f)(x) \leq \liminf_{n \to \infty}(\bar{R}_{\alpha,H}(f*\phi_n))(x) \leq \Upsilon_H(\alpha),$$
where $(\phi_n)_{n \geq 1}$ is a sequence of approximations to the identity, consisting of probability density functions in $\cS(\bR^d)$. Hence,
$$\sup_{x \in \bR^d}(\bar{R}_{\alpha,H}f)(x) \leq \Upsilon_{H}(\alpha).$$

For the reverse inequality, we let $\phi_n(x)=(2\pi)^{-d/2}n^{d/2}\exp(-n|x|^2/2)$.
By Lemma \ref{max-principle-step2},
$$(\bar{R}_{\alpha,H}(f*\phi_n))(0)=\alpha_H \int_{\bR^d} e^{-|\xi|^2/(2n)}\int_{\bR_{+}^2}|r-s|^{2H-2} e^{-(\alpha+2{\rm Re}\Psi(\xi))(r+s)} dr ds \mu(d\xi),$$
and therefore, by applying the monotone convergence theorem,
\begin{equation}
\label{limit-R-n-0}
\lim_{n \to \infty}(\bar{R}_{\alpha,H}(f*\phi_n))(0)=\Upsilon_{H}(\alpha).
\end{equation}

Using (\ref{calcul-P-st}) and the symmetry of the function $\phi_n$, we obtain:
$$(\bar{P}_{r+s}(f*\phi_n))(0)=\int_{\bR^d}(\bar{P}_{r+s}f)(x)\phi_n(x) dx.$$
Therefore,
\begin{equation}
\label{calcul-R-n-0}
(\bar{R}_{\alpha,H}(f*\phi_n))(0)=\int_{\bR^d}(\bar{R}_{\alpha,H}f)(x)\phi_n(x) dx \leq \sup_{x \in \bR^d} (\bar{R}_{\alpha,H}f)(x).
\end{equation}
From (\ref{limit-R-n-0}) and (\ref{calcul-R-n-0}), we infer that
$\Upsilon_{H}(\alpha) \leq \sup_{x \in \bR^d} (\bar{R}_{\alpha,H}f)(x)$.

The last assertion follows by taking $\phi_n$ with the support in the ball of radius $1/n$ and center $0$. $\Box$

\begin{lemma}
\label{max-principle-step4}
If $f$ satisfies (\ref{cond-on-f}), then
for any $\alpha>0$,
$$(\bar{R}_{\alpha,H}f)(0)=\Upsilon_{H}(\alpha).$$
\end{lemma}

\noindent {\bf Proof:}
Using (\ref{harmonic-result}), we have:
$$(\bar{P}_{r+s}f)(0)=\int_{\bR^d}\int_{\bR^d}\bar{p}_r(x)\bar{p}_{s}(y)f(x-y)dx dy=\int_{\bR^d}e^{-2(r+s){\rm Re} \Psi(\xi)}g(\xi)d\xi.$$
The conclusion follows from the definitions of $(\bar R_{\alpha,H}f)(0)$ and $\Upsilon_{H}(\alpha)$.
$\Box$

\vspace{3mm}

To investigate the connection with the parabolic problem (\ref{parabolic-eq}), we let $$\Upsilon_H^{*}(\alpha)=\int_{\bR^d} \left(\frac{1}{\alpha+2{\rm Re}\Psi(\xi)} \right)^{2H}\mu(d\xi).$$
By Lemma \ref{lemmaA} (Appendix A), $\Upsilon_H^{*}(\alpha)<\infty$ for all $\alpha>0$ if and only if $\Upsilon_H^{*}(\alpha)<\infty$ for some $\alpha>0$.

The following result gives the relationship between $\Upsilon_H(\alpha)$ and $\Upsilon_H^*(\alpha)$.

\begin{lemma}
\label{ineq-Upsilon} For any $\alpha>0$,
$$c_{\alpha,H} \Upsilon_H^*(\alpha) \leq \Upsilon_{H}(\alpha) \leq b_H^2 H^{2H}
\Upsilon_{H}^{*}(\alpha),$$
where $c_{\alpha,H}=2^{-(2H+2)}[(\alpha \wedge 1)/(\alpha+3/2)]^{2H}$.
\end{lemma}

\noindent {\bf Proof:} The second inequality follows by (\ref{MMV-ineq}).
For the first inequality, we note that, since the integrand from the definition of $\Upsilon_H(\alpha)$ is non-negative, the integral $dr ds$ over $[0,\infty)^2$ can be bounded below by the integral over $[0,1]^2$. By Proposition 4.3 of \cite{BT10-SPA}, for any $t>0$ and $\lambda \geq 0$
$$\alpha_{H} \int_0^t \int_0^t |r-s|^{2H-2} e^{-\lambda (r+s)} dr ds \geq \frac{1}{4} (t^{2H} \wedge 1) \left( \frac{1}{2}\right)^{2H}\left( \frac{1}{1/2+\lambda}\right)^{2H}.$$
Applying this result for $t=1$ and $\lambda=\alpha+2{\rm Re} \Psi(\xi)$, we obtain:
$$\int_0^1 \int_0^1 |r-s|^{2H-2} e^{-(\alpha+2{\rm Re} \Psi(\xi))(r+s)} dr ds \geq \left( \frac{1}{2}\right)^{2H+2} \left(\frac{1}{1/2+\alpha+2{\rm Re} \Psi(\xi)} \right)^{2H}.$$
Hence,
$$\Upsilon_H(\alpha) \geq \left(\frac{1}{2}\right)^{2H+2} \Upsilon_H^*(\alpha+1/2) \geq c_{\alpha,H} \Upsilon_H^*(\alpha),$$
where we used Lemma \ref{lemmaA} (Appendix A) for the second inequality. $\Box$

\vspace{3mm}

Recall that by Theorem \ref{th-existence-parabolic}, condition (\ref{parabolic-cond}) is the necessary and sufficient for problem (\ref{parabolic-eq}) to have a random field solution. As a consequence of the maximum principle, we obtain the following result.

\begin{corollary}
\label{resolvent-cond-equiv}
Suppose that $f$ satisfies (\ref{cond-on-f}).
Then (\ref{parabolic-cond}) is equivalent to
\begin{equation}
\label{resolvent-cond} (\bar{R}_{\alpha,H}f)(0)<\infty \quad
\mbox{for any} \ \alpha>0.
\end{equation}
\end{corollary}

\noindent {\bf Proof:} By Theorem \ref{max-principle}, (\ref{resolvent-cond}) holds
if and only if $\Upsilon_H(\alpha)<\infty$ for any $\alpha>0$. By Lemma \ref{ineq-Upsilon}, this is equivalent to $\Upsilon_H^*(\alpha)<\infty$ for any $\alpha >0$, which in turn, is equivalent to
(\ref{parabolic-cond}) (i.e. $\Upsilon_H^*(2)<\infty$), by Lemma \ref{lemmaA} (Appendix A).  $\Box$

\subsection{Connection with the Intersection Local Time}

When the noise $W$ is white in time, the authors of \cite{FKN09} and \cite{FK10} noticed an interesting connection between the existence of a random field solution of problem (\ref{parabolic-eq}) and the existence of the occupation time
$$L_t(f)=\int_0^t f(\bar{X}_s)ds.$$


In this section, we develop a similar connection in the case of the fractional noise, by considering the ``weighted'' intersection local time:
$$L_{t,H}(f)=\alpha_H \int_0^t \int_0^t |r-s|^{2H-2}f(\bar{X}_r^1-\bar{X}_s^2)dr ds,$$
where $(\bar{X}_t^1)_{t \geq 0}$ and $(\bar{X}_t^2)_{t \geq 0}$ are two independent copies of $(\bar{X}_t)_{t \geq 0}$.


Clearly for any fixed $t>0$, $E[L_{t,H}(f)]<\infty$ is a sufficient condition for $L_{t,H}(f)<\infty$ a.s., but the negligible set depends on $t$. Our result will show that under condition (\ref{resolvent-cond}), $L_{t,H}(f) <\infty$ for all $t>0$ a.s.
To motivate this result, we consider first an example, in which we proceed to the calculation of $E[L_{t,H}(f)]$ in a particular case.

Note that $\bar{X}_r^1-\bar{X}_s^2 \stackrel{d}{=}\bar{X}_{r+s}$ for any $r,s \in [0,t]$,
and therefore,
\begin{equation}
\label{mean-Lt(f)}
E[L_{t,H}(f)]=\alpha_H \int_0^t \int_0^t |r-s|^{2H-2}E[f(\bar{X}_{r+s})]dr ds.
\end{equation}

\begin{example}{\rm (Stable processes)} {\rm Refer to Example \ref{stable-ex}.
Since $(\bar{X}_t)_{t \geq 0}$ is self-similar with exponent $1/\beta$ (see Theorem 13.5 of \cite{sato99}),
$$\bar{X}_{r+s}\stackrel{d}{=} (r+s)^{1/\beta}\bar{X}_{1}.$$

Suppose in addition that $f(x)=|x|^{-(d-\alpha)}$ for $0<\alpha<d$. Then $E[f(\bar{X}_{r+s})]=E|\bar{X}_{r+s}|^{-(d-\alpha)}=c_{\alpha,d}(r+s)^{-(d-\alpha)/\beta}$, where
$c_{\alpha,d}=E|\bar{X}_1|^{-(d-\alpha)}$. By (\ref{mean-Lt(f)}),
$$E[L_{t,H}(f)]=\alpha_H c_{\alpha,d}\int_0^t \int_0^t |r-s|^{2H-2}(r+s)^{-(d-\alpha)/\beta}dr ds.$$
One can see that $E[L_{t,H}(f)]<\infty$ for any $t>0$ if and only if
(\ref{cond-H-alpha-beta}) (or equivalently, (\ref{parabolic-cond})) holds. But by Corollary \ref{resolvent-cond-equiv},
(\ref{parabolic-cond}) is equivalent to (\ref{resolvent-cond}). }
\end{example}

The previous example shows that the existence of $L_{t,H}(f)$ is related to the
potential-theoretic condition (\ref{resolvent-cond}),  which is in turn the necessary and sufficient condition for the existence of  a random field solution to problem (\ref{parabolic-eq}) (by Theorem \ref{th-existence-parabolic} and Corollary \ref{resolvent-cond-equiv}).
We will see below that this is a general phenomenon. For this, suppose that $$\bar{X}_0^1=x_1, \quad \bar{X}_0^2=x_2,$$
and let $P_{x_i}$ be the law of $\bar{X}^i$ for $i=1,2$. Then $P_{x_1,x_2}=P_{x_1} \times P_{x_2}$ is the law of $(\bar{X}^1,\bar{X}^2)$. We denote by $E_{x_1,x_2}$ the expectation under $P_{x_1,x_2}$.

The next result shows the existence of the intersection local time $L_{t,H}(f)$ under condition (\ref{resolvent-cond}).

\begin{theorem} {\rm (Connection with the Local Time)}
Suppose that $f$ satisfies (\ref{cond-on-f}). If
(\ref{resolvent-cond}) holds, then for any $x_1,x_2 \in \bR^d$,
$$P_{x_1,x_2}(L_{t,H}(f)<\infty \quad \mbox{for all} \ t>0)=1$$
$$P_{x_1,x_2}\left(\limsup_{t \to \infty}\frac{\log L_{t,H}(f)}{t} \leq 0\right)=1.$$
\end{theorem}

\noindent {\bf Proof:} We follow the lines of the proof of Theorem 3.13 of \cite{FK10}. Since $f$ is non-negative, it follows that for any $t>0$,
\begin{equation}
\label{bound-Lt-e}
e^{-2\alpha t} L_{t,H}(f) \leq \alpha_H \int_0^{\infty} \int_0^{\infty} e^{-\alpha (r+s)}|r-s|^{2H-2}f(\bar{X}_r^1-\bar{X}_s^2)drds.
\end{equation}

Note that
$$E_{x_1,x_2}[f(\bar{X}_r^1-\bar{X}_s^2)]=\int_{\bR^d} \int_{\bR^d}f(y-z)\bar{p}_{r}(x_1-y)\bar{p}_{s}(x_2-z)dy dz=(\bar{P}_{r+s}f)(x_1-x_2).$$
Taking supremum over $t$, and expectation with respect to $P_{x_1,x_2}$ in (\ref{bound-Lt-e}), we obtain:
$$E_{x_1,x_2}[\sup_{t>0}(e^{-2\alpha t} L_{t,H}(f))] \leq (\bar{R}_{\alpha,H}f)(x_1-x_2).$$
From here, using Theorem \ref{max-principle} and condition (\ref{resolvent-cond}), we infer that:
$$\sup_{x_1,x_2 \in \bR^d}E_{x_1,x_2}[\sup_{t>0}(e^{-2\alpha t} L_{t,H}(f))]  \leq \sup_{x \in \bR^d} (\bar{R}_{\alpha,H}f)(x)=(\bar{R}_{\alpha,H}f)(0)<\infty.$$
The result follows. $\Box$

\section{The Hyperbolic Equation}
\label{hyperbolic-section}

In this section we consider the hyperbolic problem
(\ref{hyperbolic-eq}). Throughout this section, we assume that $X$ is symmetric, i.e.
\begin{equation}
\label{imaginary-zero}
{\rm Im} \Psi(\xi)=0, \quad \mbox{for all} \ \xi \in \bR^d.
\end{equation}
Since ${\rm Re} \Psi(\xi)=\Psi (\xi)$, we use the notation $\Psi(\xi)$
to simplify the writing.

To define the weak solution, we cannot use the same method as in the parabolic case, since in general, we may not be able to identify the fundamental solution $G$ of $\partial_{tt}u-\cL u=0$. Note that in some particular cases, we are able to identify $G$ (see Remark \ref{stable-rem-hyp} below).

 To circumvent this difficulty, we use the method of \cite{FKN09}, whose salient features we recall briefly below. Consider first the deterministic equation:
\begin{equation}
\label{deterministic-wave}
\frac{\partial^2 u}{\partial t^2}(t,x)=\cL u(t,x)+F(t,x),
\end{equation}
with zero initial conditions, where $F$ is a smooth function.
By taking formally the Fourier transform in the $x$ variable, and using the fact that $\cF {\cL}=-\Psi$, we obtain that $\cF u$ satisfies the following equation:
\begin{equation}
\label{deterministic-wave-Fourier}
\frac{\partial^2 (\cF u)}{\partial t^2}(t,\xi)=-\Psi(\xi)\cF u(t,\xi)+\cF F(t,\xi),
\end{equation}
with zero initial conditions. Equation (\ref{deterministic-wave-Fourier}) can be solved using Duhamel's principle. We obtain:
$$\cF u(t,\xi)=\frac{1}{\sqrt{\Psi(\xi)}} \int_0^t \sin (\sqrt{\Psi(\xi)}(t-s)) \cF F(s,\xi)ds.$$
We apply formally the Fourier inversion formula. Multiplying by $\varphi \in \cS(\bR^d)$, and integrating $dx$, we arrive to the following (formal) definition of a weak solution of (\ref{deterministic-wave}):
\begin{equation}
\label{formal-solution-wave}
u(t,\varphi)=\frac{1}{(2\pi)^{d}} \int_{0}^{t}\int_{\bR^d}\frac{\sin (\sqrt{\Psi(\xi)}(t-s))}{\sqrt{\Psi(\xi)}}\overline{\cF \varphi(\xi)} \cF F(s,\xi) d\xi ds.
\end{equation}

If instead of the smooth function $F$ we consider the random noise
$\dot W$,
 the integral above is replaced by a stochastic integral $\cF W(ds,d\xi)$, where $\cF W$ is a Gaussian process which we define below.

As in \cite{balan-tudor08}, we let $\cP(\bR^d)$ be the completion of $\cS(\bR^d)$ with respect to the inner product:
\begin{eqnarray*}
\langle \varphi_1, \varphi_2 \rangle_{\cP(\bR^d)} &:=&\int_{\bR^d} \int_{\bR^d} \varphi_1(x)\overline{\varphi_2(y)}f(x-y)dx dy \\
&=& \int_{\bR^d} \cF \varphi_1(\xi) \overline{\cF \varphi_2(\xi)}\mu(d\xi).
\end{eqnarray*}

Let $\widehat{\cP(\bR^d)}$ be the completion of $\cS(\bR^d)$  with respect to the inner product:
\begin{equation}
\label{def-product-hat-P}
\langle \psi_1,\psi_2 \rangle_{\widehat{\cP(\bR^d)}}:=
\langle \cF \psi_1, \cF \psi_2 \rangle_{\cP(\bR^d)}.
\end{equation}

Note that if the noise $W$ is white in space, then by Plancherel
theorem, $\langle \psi_1,\psi_2
\rangle_{\widehat{\cP(\bR^d)}}=(2\pi)^d\langle \psi_1,\psi_2
\rangle_{L^2(\bR^d)}$ and
$\widehat{\cP(\bR^d)}=\cP(\bR^d)=L^2(\bR^d)$.

The following lemma gives a more direct way of calculating $\langle \psi_1,\psi_2 \rangle_{\widehat{\cP(\bR^d)}}$.

\begin{lemma}
\label{direct-calcul-norm-hatHP}
For any $\psi_1,\psi_2 \in \cS(\bR^d)$,
$$\langle \psi_1,\psi_2 \rangle_{\widehat{\cP(\bR^d)}}=(2\pi)^{2d} \int_{\bR^d} \psi_1(\xi)\overline{\psi_2(\xi)}\mu(d\xi).$$
\end{lemma}

\noindent {\bf Proof:} Note that for any $\varphi \in L^1(\bR^d)$, $\overline{\cF \varphi(\xi)}=\cF^{-1}\overline{\varphi}(\xi)$, where
$$\cF^{-1} \varphi(\xi):=\int_{\bR^d} e^{i \xi \cdot x}\varphi(x)dx, \quad \forall \xi \in \bR^d.$$

We denote $\varphi_i:=\overline{\cF \psi_i} \in \cS(\bR^d)$ for $i=1,2$. We obtain:
\begin{eqnarray*}
\langle \psi_1,\psi_2 \rangle_{\widehat{\cP(\bR^d)}} &=&
\int_{\bR^d} \int_{\bR^d} \overline{\varphi_1(x)}\varphi_2(y)f(x-y)dx dy \\
&=& 
\int_{\bR^d} \overline{\cF \varphi_1(\xi)}\cF \varphi_2(\xi) \mu(d\xi) \\
&=& 
\int_{\bR^d} \cF^{-1} (\cF \psi_1)(\xi) \overline{\cF^{-1}(\cF \psi_2)(\xi)} \mu(d\xi).
\end{eqnarray*}
By the Fourier inversion theorem, $\cF^{-1}(\cF \psi_i)=(2\pi)^d \psi_i$ for $i=1,2$. The result follows. $\Box$

We endow the space $\cE$ with the inner product:
$$\langle h_1,h_2 \rangle_{\widehat{\cH \cP}}:=
\langle \cF h_1, \cF h_2 \rangle_{\cH \cP},$$
where $\cF$ denotes the Fourier transform in the $x$ variable.

Note that by Lemma \ref{direct-calcul-norm-hatHP}, for any $h_1, h_2 \in \cE$,
$$\langle h_1,h_2 \rangle_{\widehat{\cH \cP}}= \alpha_H (2\pi)^{2d} \int_{\bR^d} \int_{0}^{\infty} \int_{0}^{\infty}|r-s|^{2H-2} h_1(r,\xi) \overline{h_2(s,\xi)}  dr ds\mu(d\xi).$$



For any $h \in \cE$, we define:
$$\cF W(h):=
W(\cF h).$$
By the isometry property of $W$, for any $h \in \cE$,
$$E|\cF W(h)|^2=E|W(\cF h)|^2= \| \cF h \|_{\cH \cP}^2=\|h\|_{\widehat{\cH \cP}}^2.$$

Let $\widehat{\cH \cP}$ be the completion of $\cE$ with respect to the inner product $\langle \cdot, \cdot \rangle_{\widehat{\cH \cP}}$. The map $\cE \ni h \mapsto \cF W(h) \in L^2(\Omega)$ is an isometry which can be extended to $\widehat{\cH \cP}$. We denote this extension by:
$$h \mapsto \int_0^{\infty}\int_{\bR^d}h(t,\xi)\cF W (dt,d\xi)=:\cF W(h).$$

This gives the rigorous construction of the isonormal Gaussian process $\cF W=\{\cF W(h); h \in \widehat{\cH \cP}\}$ which was mentioned formally above.

The following result gives a criterion for a function $h$ to be in $\widehat{\cH \cP}$.

\begin{lemma}
\label{criterion-phi-in-hatHP} Let $h: \bR_{+} \times \bR^d \to \bC$
be a deterministic function such that $h(t,\cdot)=0$ if $t>T$.
Suppose that $h$
satisfies the following conditions: \\
(i) $h(t, \cdot) \in L^2(\bR^d)$ for all $t \in [0,T]$;\\
(ii) $\cF h \in \cH \cP$, where $\cF h$ denotes the Fourier transform in the $x$ variable. \\
Then $h \in \widehat{\cH \cP}$ and
$$\|h\|_{\widehat{\cH \cP}}^2= \alpha_H (2\pi)^{2d} \int_{\bR^d} \int_{0}^{T} \int_{0}^{T}|r-s|^{2H-2} h(r,\xi) \overline{h(s,\xi)}dr ds \mu(d\xi).$$
In particular, the stochastic integral of $h$ with respect to the noise $\cF W$ is well-defined.
\end{lemma}

\noindent {\bf Proof:} Let $g=\cF h$. By the Fourier inversion formula on $L^2(\bR^d)$, $h(s,\xi)=(2\pi)^{-d}\overline{\cF(s,\xi)}$, and hence,
\begin{equation}
\label{F-inversion-g}
\cF g(s,\xi)=(2\pi)^d \overline{h(s,\xi)}.
\end{equation}
Since $g \in \cH \cP$, there exists a sequence $(g_n)_{n \geq 1}$ of the form $g_n(s,x)=\phi_n(s)\gamma_n(x)$, where $\phi_n$ is a linear combination of indicator functions $1_{[0,a]}, a \in [0,T]$ and $\gamma_n \in \cS(\bR^d)$, such that $\|g_n-g\|_{\cH \cP} \to 0$ (see \cite{BT10-SPA}).

Let $\psi_n:=(2\pi)^{-d}\overline{\cF \gamma_n} \in \cS(\bR^d)$ and $h_n(s,\xi)=\phi_n(s)\psi_n(\xi)$.
Then 
\begin{equation}
\label{F-inversion-gn}\cF g_n(s,\xi)=\phi_n(s) \cF \gamma_n(\xi)=(2\pi)^d \phi_n(s) \overline{\psi_n(\xi)}=(2\pi)^d \overline{h_n(s,\xi)}.
\end{equation}

Using (\ref{F-inversion-g}) and (\ref{F-inversion-gn}), we obtain that
\begin{eqnarray*}
\lefteqn{\alpha_H (2\pi)^{2d} \int_{\bR^d} \int_0^T \int_0^T |r-s|^{2H-2} (h_n-h)(r,\xi) \overline{ (h_n-h)(s,\xi)}dr ds \mu(d\xi) } \\
& & = \alpha_H \int_{\bR^d} \int_0^T \int_0^T |r-s|^{2H-2} \overline{(\cF g_n-\cF g)(r,\xi)} (\cF g_n-\cF g)(s,\xi) dr ds \mu(d\xi)\\
& & =\|g_n-g\|_{\cH \cP}^2 \to 0.
\end{eqnarray*}
The conclusion follows. $\Box$

\vspace{3mm}

We now return to equation (\ref{hyperbolic-eq}). By analogy with (\ref{formal-solution-wave}), we say that the process $\{u(t,\varphi); t \geq 0,\varphi \in \cS(\bR^d)\}$ defined by:
$$u(t,\varphi)=\frac{1}{(2\pi)^{d}} \int_{0}^{t} \int_{\bR^d}\frac{\sin (\sqrt{\Psi(\xi)}(t-s))}{\sqrt{\Psi(\xi)}}\overline{\cF \varphi(\xi)} \cF W(ds,d\xi),$$
is a {\em weak solution} of (\ref{hyperbolic-eq}).
The stochastic integral above is well-defined if and only if the integrand
$$(s,\xi) \mapsto h_{t,\varphi}(s,\xi)=\frac{1}{(2\pi)^d} 1_{[0,t]}(s)
\frac{\sin (\sqrt{\Psi(\xi)}(t-s))}{\sqrt{\Psi(\xi)}}\overline{\cF \varphi(\xi)}
\quad \mbox{belongs to $\widehat{\cH \cP}$}.$$

To check that $h_{t,\varphi} \in
\widehat{\cH \cP}$, it suffices to show that $h_{t,\varphi}$ satisfies conditions (i) and (ii) of Lemma \ref{criterion-phi-in-hatHP}. Condition (i) holds since $|\sin x| \leq |x|$ for any $x$. For (ii), we have to show that $g_{t,\varphi}:=\cF h_{t,\varphi} \in \cH \cP$. For this we apply Theorem 2.1 of \cite{BT10-SPA}. Note that the function $(s,\xi) \mapsto \cF g_{t,\varphi}(s,\xi)=(2\pi)^d \overline{h_{t,\varphi}(s,\xi)}$ satisfies conditions (i)-(iii) of this theorem. So, if suffices to show that:
$$I_{t,\varphi}:=\alpha_H (2\pi)^{2d}\int_{\bR^d} \int_0^{\infty}
\int_0^{\infty}|r-s|^{2H-2} \overline{h_{t,\varphi}(r,\xi)}
h_{t,\varphi}(s,\xi) dr ds \mu(d\xi)<\infty.$$

Note that
$$I_{t,\varphi}=
\int_{\bR^d}N_t(\xi)|\cF \varphi(\xi)|^2 \mu(d\xi),$$
where
$$N_t(\xi)=\frac{\alpha_H}{\Psi(\xi)} \int_0^t \int_0^t \sin (r\sqrt{\Psi(\xi)})\sin (s\sqrt{\Psi(\xi)})|r-s|^{2H-2} dr ds.$$

Using the fact that $|\sin x| \leq |x|$ for any $x$, it follows
that $$N_t (\xi) \leq \alpha_H \int_0^t \int_0^t rs |r-s|^{2H-2}drds
\leq t^{2H+2},$$ and hence $I_{t,\varphi} \leq
 t^{2H+2} \int_{\bR^d}|\cF \varphi(\xi)|^2 \mu(d\xi)<\infty$. This
proves that $u(t,\varphi)$ is a well-defined random variable in
$L^2( \Omega)$, for any $t>0$ and $\varphi \in \cS(\bR^d)$. Moreover,
$$E|u(t,\varphi)|^2=E|\cF W(h_{t,\varphi})|^2=\|h_{t,\varphi}\|_{\widehat{\cH \cP}}^2=I_{t,\varphi}.$$

To define the random field solution of (\ref{hyperbolic-eq}), we
proceed as in the case of the parabolic equation. We define the
norms:
\begin{eqnarray*}
\|\varphi\|_t^2 &:=& E|u(t,\varphi)|^2 =\int_{\bR^d} N_t(\xi) |\cF
\varphi(\xi)|^2 \mu(d\xi)\\
|\|\varphi \||_{t}^2 &:=&  \cE(t;\varphi)=\int_{\bR^d}
\left(\frac{1}{1/t+\Psi(\xi)} \right)^{H+1/2}|\cF
\varphi(\xi)|^2\mu(d\xi).
\end{eqnarray*}
Let $M_t$ and $\cZ_t$ be the completions of $\cS(\bR^d)$ with
respect to the norms $\|\cdot\|_t$, respectively $|\|\cdot \||_t$.
Let $M=\cap_{t>0}M_t$ and $\cZ=\cap_{t>0}\cZ_t$. By Lemma \ref{lemmaA} (Appendix A),
$\cZ_t=\cZ_s=\cZ$ for any $s,t>0$.

 We say that equation (\ref{hyperbolic-eq}) has a
{\em random field solution} if $\delta_x \in M$ for any $x \in
\bR^d$. In this case, the random field solution is defined by
$\{u(t,x)=u(t,\delta_x); t \geq 0,x \in \bR^d\}$.

The following result gives some upper and lower bounds for
$N_t(\xi)$. For the upper bound, we use an argument similar to
Proposition 3.7 of \cite{BT10-SPA}. For the lower bound, we use a
new argument.

\begin{theorem}
\label{bounds-for-Ntxi-hyp}
For any $t>0$ and $\xi \in \bR^d$,
$$D_H^{(2)} t \left(\frac{1}{1/t^2+\Psi(\xi)} \right)^{H+1/2} \leq
N_t(\xi) \leq D_H^{(1)} t \left(\frac{1}{1/t^2+\Psi(\xi)}
\right)^{H+1/2},$$ where $D_H^{(1)}$ and $D_H^{(2)}$ are some positive constants
depending only on $H$.
\end{theorem}

\noindent {\bf Proof:} We first prove the upper bound. Suppose that
$t^2 \Psi(\xi) \leq 1$. Using (\ref{MMV-ineq}), the fact that
$\|\varphi\|_{L^{1/H}(0,t)}^2 \leq t^{2H-1}\|\varphi\|_{L^2(0,t)}^2$
and $|\sin x| \leq x$ for all $x>0$, we obtain:
\begin{eqnarray*}
N_t(\xi) & \leq & b_H^2 t^{2H-1}\frac{1}{\Psi(\xi)} \int_0^t \sin^2
(r
\sqrt{\Psi(\xi)}) dr \leq b_H^2 t^{2H-1} \int_0^t r^2 dr \\
&=& \frac{1}{3} b_H^2 t^{2H+2} \leq \frac{1}{3} b_H^2 2^{H+1/2}t
\left(\frac{1}{1/t^2+\Psi(\xi)} \right)^{H+1/2},
\end{eqnarray*}
where for the last inequality, we used the fact that $\frac{t^2}{2}
\leq \frac{1}{t^{-2}+\Psi(\xi)}$ if $t^2 \Psi(\xi) \leq 1$.

Suppose next that $t^2 \Psi(\xi) \geq 1$. We denote $$T=t \sqrt{\Psi(\xi)}.$$ Using the change of
variable $r'=r\sqrt{\Psi(\xi)}$ and $s'=s \sqrt{\Psi(\xi)}$, we
obtain:
$$N_t(\xi)=\frac{1}{\Psi(\xi)^{H+1}} \| \sin (\cdot)\|_{\cH(0,T)}^2.$$

We now use Lemma B.1 of \cite{BT10-SPA} for expressing the
$\cH(0,T)$-norm of the sinus function in the spectral domain. We
obtain that:
\begin{equation}
\label{Nt-wave-spectral} N_t(\xi)=\frac{c_H}{\Psi(\xi)^{H+1}}
\int_{\bR} \frac{|\tau|^{-(2H-1)}}{(\tau^2-1)^2}
[f_{T}^{2}(\tau)+g_T^2 (\tau)]d\tau,
\end{equation}
where
$$f_T(\tau)=\sin(\tau T)-\tau \sin T \quad \mbox{and}
\quad g_T(\tau)=\cos(\tau T)- \cos T.$$

Letting $\varphi(x)=\sin x$, we have: (see the proof of Lemma B.1 of \cite{BT10-SPA}) $$|\cF_{0,T}\varphi(\tau)|^2=\frac{1}{(\tau^2-1)^2}[f_{T}^2(\tau)+g_T^2(\tau)].$$

We split the integral in (\ref{Nt-wave-spectral}) into the regions $|\tau| \leq 1/2$ and
$|\tau| \geq 1/2$, and denote the two integrals by $N_t^{(1)}(\xi)$
and $N_t^{(2)}(\xi)$. Using the same argument as in the proof of
Proposition 3.7 of \cite{BT10-SPA}, we get:
\begin{equation}
\label{estimate-Nt1} N_t^{(1)}(\xi) \leq C
\frac{c_H}{\Psi(\xi)^{H+1}} \cdot \frac{2^{2H-2}}{1-H}\leq C
\frac{c_H }{1-H} 2^{2H-2} t \left(
\frac{2}{1/t^2+\Psi(\xi)}\right)^{H+1/2},
\end{equation}
where $C=11.11$, and for the second inequality we used the fact that
$\frac{1}{\Psi(\xi)^{1/2}} \leq t$ and $\frac{1}{\Psi(\xi)} \leq
\frac{2}{t^{-2}+\Psi(\xi)}$ if $t^2 \Psi(\xi) \geq 1$.

On the other hand,
$$N_t^{(2)}(\xi) \leq c_H 2^{2H-1} \frac{1}{\Psi(\xi)^{H+1}}
I(T),$$ where
\begin{equation}
\label{Planch-id-wave} I(T):=\int_{\bR}
\frac{f_T^2(\tau)+g_T^2(\tau)}{(\tau^2-1)^2}d\tau
=2\pi \int_0^T \sin^2 x dx=\pi
T \left[1-\frac{\sin(2 T)}{2 T} \right].
\end{equation}
by Plancherel's
theorem. 
This yields the estimate $I(T) \leq 2\pi T$. We
obtain:
\begin{equation}
\label{estimate-Nt2} N_t^{(2)}(\xi) \leq c_H 2^{2H-1} 2\pi t
\left(\frac{1}{\Psi(\xi)}\right)^{H+1/2} \leq c_H 2^{2H-1} 2 \pi t
\left( \frac{2}{1/t^2+\Psi(\xi)} \right)^{H+1/2},
\end{equation}
where for the second inequality we used the fact that $t^2 \Psi(\xi)
\geq 1$.

Combining (\ref{estimate-Nt1}) and (\ref{estimate-Nt2}), we conclude
that:
$$N_t(\xi) \leq C \frac{c_H}{1-H} 2^{3H-1/2} t\left( \frac{1}{1/t^2+\Psi(\xi)}
\right)^{H+1/2}.$$ The upper bound follows, letting
$$D_H^{(1)}=\max \left\{ \frac{1}{3} b_H^2 2^{H+1/2}, C \frac{c_H}{1-H} 2^{3H-1/2} \right\}.$$

We now treat the lower bound. Suppose first that $t^2 \Psi(\xi) \leq
1$. Using the fact that $\sin x \geq x \sin 1$ for all $x \in
[0,1]$, we obtain:
\begin{eqnarray*}
N_t(\xi) & \geq & \alpha_H \sin^2 1 \int_0^t \int_0^t rs
|r-s|^{2H-2}dr ds= \alpha_H \sin^2 1 \frac{\beta(2,2H-1)}{H+1}
t^{2H+2} \\
& \geq & \alpha_H \sin^2 1 \frac{\beta(2,2H-1)}{H+1} t \left(
\frac{1}{1/t^2+\Psi(\xi)}\right)^{H+1/2},
\end{eqnarray*}
where $\beta$ denotes the Beta function and we used the fact that
$t^2 \geq \frac{1}{t^{-2}+\Psi(\xi)}$.

Suppose next that $T^2=t^2 \Psi(\xi) \geq 1$.
Let $\rho>1$ be a constant
which will be specified below. We use (\ref{Nt-wave-spectral}).
Since the integrand is non-negative, $N_t(\xi)$ is bounded below by
the integral over the region $|\tau| \leq \rho$. In that region,
$|\tau|^{-(2H-1)} \geq \rho^{-(2H-1)}$, and hence
\begin{equation}
\label{LB-Nt-wave-step1} N_t(\xi) \geq \frac{c_H
\rho^{-(2H-1)}}{\Psi(\xi)^{H+1}} \left(
I(T)-\int_{|\tau| \geq \rho}
\frac{f_T^2(\tau)+g_T^2(\tau)}{(\tau^2-1)^2}d\tau\right).
\end{equation}

Using (\ref{Planch-id-wave}) and the inequality $1-(\sin x)/x \geq
1/2$ for any $x \geq 2$, we get the lower bound:
\begin{equation}
\label{LB-Nt-wave-step2} I(T) \geq \frac{\pi }{2} T
, \quad \mbox{since} \quad T \geq 1.
\end{equation}

To find an upper bound for the second integral in the right-hand
side of (\ref{LB-Nt-wave-step1}), we use the fact that:
$$f_T^2(\tau)+g_T^2(\tau) \leq 2T
(1+|\tau|)^2, \quad \forall \tau \in \bR.$$

(To see this, note that $|f_T(\tau)| \leq 1+|\tau|$ and
$|f_T(\tau)| \leq 2 T|\tau|$, since $|\sin x| \leq |x|$. Hence, $f_T^2(\tau) \leq 2 T |\tau|(1+|\tau|)$.
Similarly, $|g_T(\tau)| \leq 2$ and $|g_T(\tau)| \leq
T(1+|\tau|)$, since $|1-\cos x| \leq |x|$. Hence,
$g_T^2(\tau) \leq 2T (1+|\tau|)$.)

It follows that:
\begin{equation}
\label{LB-Nt-wave-step3} \int_{|\tau| \geq \rho}
\frac{f_T^2(\tau)+g_T^2(\tau)}{(\tau^2-1)^2}d\tau
\leq  C_{\rho}T,
\end{equation}
 where $C_{\rho}=2
\int_{|\tau| \geq \rho} \frac{(1+|\tau|)^2}{(\tau^2-1)^2}d \tau$.
Using (\ref{LB-Nt-wave-step1}), (\ref{LB-Nt-wave-step2}) and
(\ref{LB-Nt-wave-step3}), we obtain that:
$$N_t(\xi) \geq \frac{c_H
\rho^{-(2H-1)}}{\Psi(\xi)^{H+1}}  \left(\frac{\pi}{2}-C_{\rho}
\right) t \sqrt{\Psi(\xi)}.$$ Choose $\rho$ large enough such that
$C_{\rho}<\pi/2$, e.g. $\rho=4$, for which $C_{\rho}<4/3$. Using the
fact that $\frac{1}{\Psi(\xi)} \geq \frac{1}{t^{-2}+\Psi(\xi)}$, we
get
$$N_t(\xi) \geq c_H
4^{-(2H-1)} \left(\frac{\pi}{2}-\frac{4}{3} \right) t
\left(\frac{1}{1/t^2+\Psi(\xi)} \right)^{H+1/2}.$$ The lower bound
follows, letting
$$D_H^{(2)}=\min \left\{\alpha_H \sin^2 1 \frac{\beta(2,2H-1)}{H+1}, c_H
4^{-(2H-1)} \left(\frac{\pi}{2}-\frac{4}{3} \right) \right\}.$$
 $\Box$

A consequence of the previous result is that the norms $\|\cdot
\|_t$ and $|\| \cdot\||_{t^2}$ are equivalent, for any $t>0$.

\begin{corollary}
For any $t>0,\varphi \in \cS(\bR^d)$,
$$d_H t\cE(t^2;\varphi) \leq E|u(t,\varphi)|^2 \leq D_H t\cE(t^2;\varphi).$$
Hence $M_{t}=\cZ_{t^2}$ for any $t>0$, and $M=\cZ$.
\end{corollary}

Below is the main result of this section.

\begin{theorem} {\rm (Existence of Solution in the
Hyperbolic Case)} Assume that (\ref{imaginary-zero}) holds. Then (\ref{hyperbolic-eq}) has a random field
solution if and only if (\ref{cond-hyp-eq}) holds.
\end{theorem}

\noindent {\bf Proof:} As in the proof of Theorem \ref{nec-suf-delta-in-Z1},
one can show that (\ref{cond-hyp-eq}) is a necessary and sufficient
condition for $\delta_x \in \cZ_1=M$ for all $x \in \bR^d$. We omit the details. $\Box$.

\begin{example}
{\rm (Stable processes)}
{\rm As in Example \ref{stable-ex}, let $\cL=-(-\Delta)^{\beta/2}$ for $\beta \in (0,2]$. Then $\Psi(\xi)=c_{\beta}|\xi|^{\beta}$ and (\ref{cond-hyp-eq}) becomes:
\begin{equation}
\label{cond-stable-proc-hyp}
\int_{\bR^d} \left(\frac{1}{1+|\xi|^{\beta}} \right)^{H+1/2}\mu(d\xi)<\infty.
\end{equation}

\noindent When $f(x)=c_{\alpha,d}|x|^{-(d-\alpha)}$ with $0<\alpha<d$, 
(\ref{cond-stable-proc-hyp}) is equivalent to
$$\left(H+\frac{1}{2}\right) \beta>d-\alpha,$$
whereas for $f(x)=\prod_{i=1}^{d}(\alpha_{H_i}|x_i|^{2H_i-2})$, 
(\ref{cond-stable-proc-hyp}) is equivalent to
$$\left(H+\frac{1}{2}\right) \beta >d-\sum_{i=1}^{d}(2H_i-1).$$
}
\end{example}

\begin{remark} \label{stable-rem-hyp}
{\rm (Fractional Powers of the Laplacian)}
{\rm As in Remark \ref{stable-rem}, we can consider the case $\cL=-(-\Delta)^{\beta/2}$ for arbitrary $\beta>0$. Note that the fundamental solution $G$ of $\partial_{tt} u- \cL u=0$ exists and satisfies:
$$\cF G (t,\xi)=\frac{\sin(t|\xi|^{\beta/2})}{|\xi|^{\beta/2}}.$$
(see p.11 of \cite{dalang-sanzsole05}).
Using Theorem 2.1 of \cite{BT10-SPA} and estimates similar to those given by Theorem \ref{bounds-for-Ntxi-hyp} above, one can show that a random field solution of (\ref{hyperbolic-eq}) (in the sense of \cite{BT10-SPA}) exists if and only if (\ref{cond-stable-proc-hyp}) holds.
}
\end{remark}

\appendix

\section{Some elementary inequalities}

\begin{lemma}
\label{lemmaA} For any $\alpha,\beta>0$,
$$c_1(\alpha,\beta)^{2H}\Upsilon_{H}^*(\beta) \leq \Upsilon_H^*(\alpha) \leq c_2(\alpha,\beta)^{2H}\Upsilon_H^*(\beta),$$
where $c_1(\alpha,\beta)=(\beta \wedge 1)/(\alpha+1)$ and $c_2(\alpha,\beta)=(\beta+1)[(1/\alpha) \vee 1]$. 
\end{lemma}

\noindent {\bf Proof:} We denote by $\Upsilon_{H,1}^{*}(\alpha)$ and $\Upsilon_{H,2}^{*}(\alpha)$ the integrals over the regions $\{2{\rm Re} \Psi(\xi) \leq 1\}$, respectively  $\{2{\rm Re} \Psi(\xi) \geq 1\}$. Using the inequality
$(\alpha+1)^{-1} \leq [1+2{\rm Re} \Psi(\xi)]^{-1} \leq \alpha^{-1}$ if
$2{\rm Re} \Psi(\xi) \leq 1$,
we obtain that:
$$\left(\frac{1}{\alpha+1}\right)^{2H} \int_{2{\rm Re} \Psi(\xi) \leq 1}\mu(d\xi) \leq \Upsilon_{H,1}^{*}(\alpha) \leq \left(\frac{1}{\alpha}\right)^{2H} \int_{2{\rm Re} \Psi(\xi) \leq 1}\mu(d\xi).$$
Combining this with the similar inequality for $\Upsilon_{H,1}^{*}(\beta)$, we get:
\begin{equation}
\label{elem-ineq-upsilon1}
\left(\frac{\beta}{\alpha+1}\right)^{2H} \Upsilon_{H,1}^{*}(\beta) \leq \Upsilon_{H,1}^{*}(\alpha) \leq \left(\frac{\beta+1}{\alpha}\right)^{2H} \Upsilon_{H,1}^{*}(\beta).
\end{equation}

A similar argument works for $\Upsilon_{H,2}^{*}(\alpha)$. We obtain:
\begin{equation}
\label{elem-ineq-upsilon2}
\left(\frac{1}{\alpha+1}\right)^{2H} \Upsilon_{H,2}^{*}(\beta) \leq \Upsilon_{H,2}^{*}(\alpha)
\leq (\beta+1)^{2H} \Upsilon_{H,2}^{*}(\beta).
\end{equation}
The result follows by taking the sum of (\ref{elem-ineq-upsilon1}) and (\ref{elem-ineq-upsilon2}).
$\Box$

\vspace{3mm}

We recall the definitions of the functionals $\cE(t;\varphi)$ introduced in Section \ref{parabolic-section}, respectively Section \ref{hyperbolic-section}. To make a distinction between these functionals in the two cases, we use the index ``par'' for parabolic, and ``hyp'' for hyperbolic:
\begin{eqnarray*}
\cE_{{\rm par}}(t;\varphi) &=& \int_{\bR^d} \left(\frac{1}{1/t+ {\rm Re} \Psi(\xi)} \right)^{2H}|\cF \varphi(\xi)|^2 \mu(d\xi)\\
\cE_{{\rm hyp}}(t;\varphi) &=& \int_{\bR^d} \left(\frac{1}{1/t+ {\rm Re }\Psi(\xi)} \right)^{H+1/2}|\cF \varphi(\xi)|^2 \mu(d\xi).
\end{eqnarray*}

\begin{lemma}
\label{ineq-cE}
For any $s>0,t>0$ and $\varphi \in \cS(\bR^d)$,
$$c_1(s,t)^{2H} \cE_{{\rm par}}(s;\varphi) \leq \cE_{{\rm par}}(t; \varphi) \leq c_2(s,t)^{2H} \cE_{{\rm par}}(s; \varphi)$$
$$c_1(s,t)^{H+1/2} \cE_{{\rm hyp}}(s;\varphi) \leq \cE_{{\rm hyp}}(t; \varphi) \leq c_2(s,t)^{H+1/2} \cE_{{\rm hyp}}(s; \varphi),$$
where $c_1(s,t)=(s^{-1} \wedge 1)/(t^{-1}+1)$ and $c_2(t)=(s^{-1}+1)(t \vee 1)$.
\end{lemma}

\noindent {\bf Proof:} The argument is similar to the proof of Lemma \ref{lemmaA}. $\Box$

\section{The $\cH(0,T)$-norm of the exponential}

The next result gives the expression of the $\cH(0,T)$-norm of the complex-valued exponential function in the spectral domain, which is needed in the proof of Theorem \ref{bounds-Nt-xi}.

\begin{lemma}
\label{H-norm-exp} Let $\varphi(x)=e^{-x(a+ib)}$ for $x \in (0,T)$, where $a,b \in \bR$. Then
$$\|\varphi\|_{\cH(0,T)}^2=c_H \int_{\bR}\frac{\sin^2[(\tau+b)T]+ \{e^{-at}-\cos[(\tau+b)T]\}^2}{a^2+(\tau+b)^2}
|\tau|^{-(2H-1)}d\tau,$$
where $c_H=\Gamma(2H+1) \sin(\pi H)/(2H)$.
\end{lemma}

\noindent {\bf Proof:} For complex-valued functions $\varphi \in L^2(0,t)$, we can apply Lemma A.1 of \cite{BT10-SPA} to ${\rm Re}\varphi$ and ${\rm Im}\varphi$ to obtain that:
\begin{equation}
\label{H-norm-spec-domain}
\|\varphi\|_{\cH(0,T)}^2=c_H \int_{\bR} |\cF_{0,T}\varphi(\tau)|^2 |\tau|^{-(2H-1)}d\tau,
\end{equation}
where  $\cF_{0,T}\varphi (\tau):=\int_0^T e^{-i \tau x}\varphi(x)dx$.

An elementary calculation shows that for $\varphi(x)=e^{-x(a+ib)}$,
$$|\cF_{0,T} \varphi(\tau)|^2=\frac{1}{a^2+(\tau+b)^2}[f_T^2 (\tau)+g_T^{2}(\tau)],$$ where
$f_T(\tau)=\sin [(\tau+b)T]$ and $g_T(\tau)=e^{-at}-\cos[(\tau+b)T]$. The result follows by (\ref{H-norm-spec-domain}). $\Box$

\section{A version of Plancherel theorem}

The following result is a version of Plancherel theorem needed for the calculation of $(\bar{P}_{r+s} (f*\phi))(x)$ in the proof of Lemma \ref{max-principle-step2}.

\begin{lemma}
\label{version-Plancherel}
For any $\varphi \in L^2(\bR^d)$, $\psi_1 \in L^2(\bR^d)$ and $\psi_2 \in L^1(\bR^d) \cap L^2(\bR^d)$,
$$\int_{\bR^d}\int_{\bR^d}\psi_1(x)\psi_2(y)\varphi(x-y)dydx=\frac{1}{(2\pi)^d}\int_{\bR^d} \cF \psi_1(\xi)
\overline{\cF \psi_2(\xi)} \ \overline{\cF \varphi(\xi)}d\xi.$$
\end{lemma}

\noindent {\bf Proof:} By Young's inequality, $\|\varphi * \psi_2\|_{2} \leq \|\varphi\|_2 \|\psi_2\|_{1}$ and $\varphi*\psi_2 \in L^2(\bR^d)$. The result follows by Plancherel theorem, since
$$\int_{\bR^d} \psi_1(x) (\varphi*\psi_2)(x)dx=\frac{1}{(2\pi)^d}\int_{\bR^d}\cF \psi_1(\xi) \overline{\cF(\varphi*\psi_2)(\xi)}d\xi.$$
$\Box$


\begin{thebibliography}{99}

\bibitem{alos-mazet-nualart01} Al\`{o}s, E., Mazet, O. and Nualart, D. (2001). Stochastic calculus with respect to Gaussian processes. {\em Ann. Probab.} {\bf 29}, 766-801.


\bibitem{balan-tudor08}  Balan, R.M. and Tudor, C. A. (2008).
The stochastic heat equation with fractional-colored noise:
existence of the solution. {\em Latin Amer. J. Probab. Math. Stat.}
{\bf 4}, 57-87.

\bibitem{BT09-JTP} Balan, R. M. and Tudor, C. A. (2010).
Stochastic heat equation with multiplicative fractional-colored
noise. {\em J. Theoret. Probab.} {\bf 23}, 834-870.


\bibitem{BT10-SPA} Balan, R. M. and Tudor, C. A. (2010). The stochastic
wave equation with fractional noise: a random field approach.
{\em Stoch. Proc. Appl.} {\bf 120}, 2468-2494.

\bibitem{carmona-coutin-montseny03} Carmona, P., Coutin, L. and Montseny, G. (2003).
Stochastic integration with respect to fractional Brownian motion. {\em Ann. Inst. H. Poincar\'e}
{\bf 39}, 27-68.

\bibitem{dalang99} Dalang, R. C. (1999). Extending martingale
measure stochastic integral with applications to spatially
homogenous s.p.d.e.'s. {\em Electr. J. Probab.} {\bf 4}, no. 6, 29 pp.

\bibitem{dalang-mueller03} Dalang, R. C. and Mueller, C. (2003).
Some non-linear s.p.d.e.'s that are second order in time. {\em
Electr. J. Probab.} {\bf 8}, 1-21.

\bibitem{dalang-quer10} Dalang, R. C. and Quer-Sardanyons, L.
(2010). Stochastic integrals for s.p.d.e.'s: a comparisson. {\em
Expositiones Mathematicae}. To appear.

\bibitem{dalang-sanzsole05} Dalang, R. C. and Sanz-Sol\'e, M. (2005). Regularity of the sample paths of a class of second order spde's. {\em J. Funct. Anal.} {\bf 227}, 304-337.

\bibitem{daprato-zabczyk92} Da Prato, G. and Zabczyk, J. (1992).
{\em Stochastic Equations in Infinite Dimensions}, Cambridge
University Press.

\bibitem{decreusefond-ustunel98} Decreusefond, L. and \"{U}st\"{u}nel, A. S. (1998).
Stochastic analysis of the fractional Brownian motion. {\em Potent. Anal.} {\bf 10}, 177-214.

\bibitem{duncan-hu-pasik00} Duncan, T. E., Hu, Y. and Pasik-Duncan, B. (2000). Stochastic calculus for fractional Brownian motion I. Theory. {\em SIAM J. Control Optim.} {\bf 38}, 582-612.

\bibitem{folland92} Folland, G. B. (1992). {\em Fourier Analysis and Its Applications}. American Mathematical Society. Providence, RI.

\bibitem{FKN09} Foondun, M., Khoshnevisan, D. and Nualart, E. (2011). A local time correspondence for stochastic partial differential equations. {\em Trans. Amer. Math. Soc.} {\bf 363}, 2481-2515.

\bibitem{FK10} Foondun, M. and Khoshnevisan, D. (2010). On the stochastic heat equation with spatially-colored random forcing. Preprint available at arXiv:1003.0348.

\bibitem{hitsuda72} Hitsuda, M. (1972).
Formula for Brownian partial derivatives. In: ``Proceedings of the
Second Japan-USSR Symposium on Probability Theory''. Vol. 2,
111-114.


\bibitem{hu-nualart09} Hu, Y. and Nualart, D. (2009). Stochastic
heat equation driven by fractional noise and local time. {\em
Probab. Theory Rel. Fields} {\bf 143}, 285-328.

\bibitem{hurst51} Hurst, H. E. (1951). Long term storage capacity in reservoirs. {\em Trans. Amer. Soc. Civil Eng.} {\bf 116}, 400-410.

\bibitem{ito44} It\^{o}, K. (1944). Stochastic integral. {\em Proc.
Imp. Acad. Tokyo} {\bf 20}, 519-524.

\bibitem{ito51} It\^{o}, K. (1951). On stochastic differential
equations. {\em Mem. Amer. Math. Soc.} {\bf 4}.

\bibitem{kabanov75} Kabanov, Yu. M. (1975). Extended stochastic
integrals (Russian). {\em Teor. Veroj. Primenen} {\bf 20},
725-737.

\bibitem{khoshnevisan-xiao09} Khoshnevisan, D. and Xiao, Y. (2009).
 Harmonic analysis of additive L\'evy processes.
 {\em Probab. Th. Rel. Fields} {\bf 145}, 459-515.

\bibitem{kolmogorov40} Kolmogorov, A. N. (1940). Wienersche Spiralen und einige andere interessante Kurven in Hilbertschen Raum. {\em C.R. (Doklady) Acad. USSR (N.S.)} {\bf 26}, 115-118.

 \bibitem{krylov99} Krylov, N. V. (1999). An analytic approach to
SPDEs. In {\em Stochastic partial differential equations: six
perspectives}, {\bf 64} {\em Math. Surveys Monogr.}, 185--242,
AMS, Providence, RI.

 \bibitem{kunita-watanabe67} Kunita, H. and Watanabe, S. (1967).
 On square-integrable martingales. {\em Nagoya J. Math.} {\bf 30},
 209-245.

\bibitem{mandelbrot-vanness68} Mandelbrot, B. B. and Van Ness, J. W. (1968).
Fractional Brownian motions, fractional noises and applications. {\em SIAM Review}
{\bf 10}, 422-437.

\bibitem{MMV01} Memin, J., Mishura, Y. and Valkeila, E. (2001).
Inequ alities for the moments of Wiener integrals with respect to
fractional Brownian motions. {\em Stat. Probab. Letters} {\bf 55},
421-430.

\bibitem{nualart06} Nualart, D. (2006). {\em The Malliavin Calculus and Related Topics}. Second Edition. Springer-Verlag, Berlin.


\bibitem{sato99} Sato, K.-I. (1999). {\em L\'evy Processes and Infinitely Divisible Distributions}.
Cambridge University Press.

\bibitem{schwartz66} Schwartz, L. (1966). {\em Th\'eorie des distributions}. Hermann, Paris.

\bibitem{skorohod75} Skorohod, A. V. (1975). A generalization of a stochastic integral. {\em Theory Probab. Appl.} {\bf 20}, 219-233.

\bibitem{stein70} Stein, E. M. (1970). {\em Singular Integrals and Differentiability Properties of
Functions}. Princeton University Press. Princeton, New Jersey.

\bibitem{walsh86} Walsh, J. B. (1986). An introduction to stochastic
partial differential equations. {\em Ecole d'Et\'{e} de
Probabilit\'{e}s de Saint-Flour XIV. Lecture Notes in Math.} {\bf
1180}, 265-439. Springer-Verlag, Berlin.


\end{thebibliography}
\end{document}